\documentclass[11pt,twoside]{article}  

\usepackage[ansinew]{inputenc}
\usepackage{amsfonts}
\usepackage{amssymb,amsmath}
\usepackage{latexsym}

\newcommand {\debeq}	{\begin{eqnarray*}}
\newcommand {\fineq}	{\end{eqnarray*}}

\newcommand	{\tendinfty}	
{\rightarrow\infty}

\newcommand	{\ZZ}{\mathbb{Z}}

\newcommand	{\NN}{\mathbb{N}}

\newcommand	{\RR}{\mathbb{R}}
\newcommand	{\PP}{\mathbb{P}}
\newcommand	{\EE}{\mathbb{E}}

\newtheorem	{thm}		{Theorem}[section]

\newtheorem     {rem}           {Remark}
\newtheorem	{prop}	[thm]{Proposition}

\begin{document}
\title{
New approaches of source-sink metapopulations decoupling the roles of demography and dispersal 
}
\author{\textsc{By Vincent Bansaye and  Amaury Lambert}
}
\date{}								
\maketitle
\noindent

\medskip

\noindent
\textsc{Vincent Bansaye\\
CMAP\\
\'Ecole Polytechnique\\
Route de Saclay\\
F-91128 Palaiseau Cedex, France}\\
\textsc{E-mail: }vincent.bansaye@polytechnique.edu\\
\textsc{URL: }http://www.cmapx.polytechnique.fr/$\sim$bansaye/
\\
\\
\noindent
\textsc{Amaury Lambert\\
UPMC Univ Paris 06\\
Laboratoire de Probabilités et Modèles Aléatoires CNRS UMR 7599\\
And\\
Collège de France\\
Center for Interdisciplinary Research in Biology CNRS UMR 7241\\
Paris, France}\\
\textsc{E-mail: }amaury.lambert@upmc.fr\\
\textsc{URL: }http://www.proba.jussieu.fr/pageperso/amaury/index.htm
\\
\\

\begin{abstract}
\noindent
Source-sink systems are metapopulations of habitat patches with different, and possibly temporally varying, habitat qualities, which are commonly used in ecology to study the fate of spatially extended natural populations. 

We propose new techniques that allow to disentangle the respective contributions of demography and dispersal to the dynamics and fate of a single species in a source-sink metapopulation. Our approach is valid for a general class of stochastic, individual-based, stepping-stone models, with density-independent demography and dispersal, provided the metapopulation is finite or else enjoys some transitivity property. 

We provide 1) a simple criterion of persistence, by studying the motion of a single random disperser until it returns to its initial position;  2) a joint characterization of the long-term growth rate and of the asymptotic occupancy frequencies of the ancestral lineage of a random survivor, by using large deviations theory. Both techniques yield formulae decoupling demography and dispersal, and can be adapted to the case of periodic or random environments, where habitat qualities are autocorrelated in space and possibly in time.


In this last case, we display examples of coupled time-averaged sinks allowing survival, as was previously known in the absence of demographic stochasticity for fully mixing \cite{JY98} and even partially mixing \cite{ERS12, Sch10} metapopulations.

\end{abstract}  	

\tableofcontents

\section{Introduction}
\subsection{Ecological background}

Stochastic models of population dynamics play a prominent role in epidemiology and in ecology \cite{KotBook}, in predicting the fate of natural populations  (persistence vs extinction, disease outbreak, species invasion, competitive coexistence,...), and in computing some related quantities of interest (extinction probability, long-term growth rate, mid-term equilibrium distribution, stable age distribution, parameter elasticities). 

Detailing these models is indispensable to understand the effect on these predictions of internal or external characteristic features, like spatial structure, age structure, intra- and interspecific interactions, or environmental change. In particular, \emph{metapopulation models} \cite{HGBook}, where the spatial structure is explicit, are used to infer the processes which have shaped contemporary range distributions, to predict migration trends or invasion fronts in response to biotic or abiotic changes, to understand the evolution of dispersal, to design protected areas and natural reserves, etc.

When the landscape is heterogeneous in terms of habitat suitability, even the mere question of predicting persistence can be a complicated task, since persistence is the result of the intricate interplay between population growth in suitable habitats, population depletion in unsuitable habitats and of how dispersal connects different habitat patches. In ecology, metapopulation models where habitat suitability is spatially heterogeneous are commonly referred to as \emph{source--sink systems} \cite{Dia96, Hol85, Pul88}. Roughly speaking, even if the definition of sources and sinks have been subject to debate \cite{Pul88, RRN06}, sources designate habitat patches where the habitat is suitable enough for the population to persist in the absence of dispersal (fundamental niche), and sinks are habitat patches where the population would become extinct in the absence of dispersal, or from which mortality during dispersal is too high to compensate growth. Spatial heterogeneity can be due to biotic environmental variables (predation risk, resource availability) or to abiotic environmental variables, which can either be constant through time (altitude or depth, latitude) or variable through time (precipitation, moisture, irradiance, pH, salinity). 

To study the persistence of a single species in a metapopulation, it is common to further assume that population dynamics are \emph{density-independent}. This assumption does certainly not hold for all natural populations, but can at least be used for populations whose persistence is guaranteed whenever their abundance is large enough to make this approximation unrealistic. It is also particularly relevant when asking about the establishment success of a new variant arising in few copies (immigrants, genetic mutants, infectives). 

The assumption of density-independence allows theoretical ecologists to make use of linear models: \emph{matrix population models} \cite{CasBook2} for deterministic dynamics, \emph{multitype branching processes} \cite{AHBook, ANBook, HJVBook, JagBook} for stochastic dynamics. These models are parsimonious in the number of parameters, and the associated mathematical theory is extremely well developed. The extinction probability has a very simple power dependence upon initial population size and composition, and under suitable assumptions, conditional on long-term survival, the geographic distribution of the population stabilizes over time, whereas its overall abundance grows exponentially with an exponent called the \emph{Malthusian growth rate}, or \emph{long-term growth rate}, or simply \emph{growth rate}. In addition, the stable geographic distribution and the long-term growth rate are solutions to a well-known spectral problem. Namely, the growth rate is the maximal eigenvalue of the mean offspring matrix (encompassing both demography and dispersal), and the stable distribution is an associated eigenvector \cite{LS02, SenBook}.

A lot of work has been dedicated to extend these results to more complicated situations, like infinite metapopulations \cite{MG01}, or, as earlier stressed, because spatial heterogeneity can itself be time-variable, to multitype branching processes in random environment \cite{AK71a, AK71b, BS09, HI96, HV03, Kap74, Tan77}. More ecologically-related work has investigated which dispersal strategies are more likely to persist in metapopulations with random environment \cite{GH02, Sch10, SLS09}, which such metapopulations are more prone to persistence \cite{BPR02}, and which introduction strategies are more successful (single large vs several small) \cite{HV03, WM85}. Specific attention has been given to coupled sinks, that is, metapopulations where each habitat patch is a (time-averaged) sink, but where populations might still persist thanks to dispersal in sparse favourable periods \cite{ERS12, JY98, RHB05, Sch10}.

\subsection{Goals and outline of the paper}

In the present paper, our aim is to develop new methods in order to disentangle the contributions of demography and dispersal to the dynamics and outcome of source-sink systems with possibly varying environment. We will be interested in criteria for global persistence and  in the computation of the long-term growth rate, and of the occupation frequencies of long-lived lineages.

One of the main problems of the spectral approach to the study of metapopulations is that the computation of eigenvalues and eigenvectors is totally opaque to biological interpretation. In particular, the respective contributions of dispersal and demography to the value of the long-term growth rate are very hard, if not impossible in general, to disentangle. As regards the question of persistence, we could ask for an alternative criterion, equivalent to, but simpler than, the positivity of this growth rate, which would avoid computing directly this eigenvalue. Similarly as in \cite{ERS12, HB06, KL10, RRN06, SLS09}, we will first seek to provide such an alternative criterion.

For example, in (st)age-structured models, it is easy to compute the net reproductive number $R_0$, which is the expected total progeny produced in the lifetime of a single individual. Then thanks to a simple renewal argument, the condition $R_0>1$ is seen to be equivalent to possible survival. More rigorously, the set of juvenile offspring of a focal juvenile ancestor forms what is called a \emph{stopping line}, for which it is known that an extended branching property holds \cite{Cha86}. This idea of the next generation-stopping line has been adapted to the spatial context in \cite{KL10, RRN06}, but remains of limited applicability. In the first part of this work, the key idea is to use as an alternative stopping line the set of descendants of a focal ancestor who are the first to return to the ancestor patch. Then by the extended branching property, the population will persist with positive probability iff the expected number, say $R$, of individuals on the stopping line is larger than $1$. 

If, as we first assume, the dispersal scheme does not depend on the state of the environment, then $R$ can be expressed separately in terms of the mean offspring numbers in each patch (and in each environmental state) and of the motion of a single random disperser. More specifically, a random disperser is a single walker on the metapopulation which follows the  dispersal stochastic scheme. We denote by $X_n$ its position at time $n$, so that $(X_n)$ is  a Markov chain with transitions given by the dispersal matrix which will be denoted  by $D$. In the case when the environment is constant, we let $m_i$ be the mean number of offspring begot in patch $i$, and we prove that
$$
R=m_1\EE\left(\prod_{n=1}^{T-1} m_{X_n}\right),
$$
where $T$ is the first time the random disperser returns to patch $1$ (assumed to be the initial patch). The population persists with positive probability iff $R>1$. This way, our formulae are seen to disentangle the effects of demography and dispersal. If all other habitat patches than patch $1$ have the same mean offspring $m$, then the last equality specializes into
$$
R=m_1\EE\left(m^{T-1}\right),
$$
where the expectation in the last display can now be seen as the probability generating function of the random variable $T-1$ evaluated at $m$. We will also compute this expectation in some special cases of interest. It is interesting to note that the formulae obtained in \cite{HB06} by a totally different method (expanding principle minors of the mean offspring matrix minus the identity matrix) feature numerous multiplicative terms also evoking closed reproductive paths.

In a second part, we will use large deviations techniques to prove that the logarithm  of the long-term growth rate $\rho$ and the asymptotic fraction $(\varphi_i)$  of time spent in each patch of the ancestral lineage of a random survivor, are given respectively by the maximum and the unique argmax of a functional $R-I$ defined on the set $\cal F$ of frequencies indexed by the metapopulation, where $R$ only depends on the reproduction/survival scheme and $I$ only depends on the dispersal scheme. Our formulae are then seen to decouple once again demography and dispersal. Namely,
$$
\log(\rho) = \sup\{R(f) - I(f): f\in {\cal F}\} = R(\varphi)-I(\varphi),
$$
where $R$ is a linear functional of frequencies only depending on the mean offspring numbers in each patch
$$
R(f):=\sum_i f_i \log(m_i),
$$
and $I$ is a (more complicated) functional which only depends on the dispersal matrix $D$
$$
I(f):=\sup\left\{\sum_{i} f_i\log(v_i/(vD)_i) : v\gg0\right\},
$$
where $v\gg0$ denotes a positive row vector, that is, $v_i>0$ for each $i$. We find that $(\varphi_i)$ never equals the stationary distribution of the single disperser, except in the case when all habitat qualities are identical (i.e., $m_i= m$ for all $i$). We compute $\rho$ and $\varphi$ in the case of a \emph{fully mixing} metapopulation, i.e., when the probability for an individual to migrate from patch $i$ to patch $j$ does not depend on $i$, a case also referred to as \emph{parent-independent migration}.

Addressing those questions is much more difficult when the model is enriched with a variable environment affecting simultaneously all habitat qualities. We can nevertheless adapt our arguments to the case when the environment is periodic or given by an ergodic sequence of random variables. We illustrate our speculations with two-patch metapopulations and a two-state environment. We make computations for fully mixing metapopulations. 

Finally, we prove that for both periodic and ergodic environments, we can find parameters for which there is possible survival in coupled sinks, a result which was previously known in the absence of demographic stochasticity for fully mixing \cite{JY98} and even partially mixing \cite{Sch10, ERS12} metapopulations.

Finally, we extend naturally our approach to a wide class of infinite metapopulations, called \emph{finitely transitive}, in the sense that they can be naturally built by connecting copies of a finite subgraph called \emph{motif}.

The paper ends with a short discussion on the uses of and possible extensions to our method (other models, local vs global persistence).

\section{Preliminaries}

\subsection{Model}

We consider a stochastic, individual-based model of spatially structured population dynamics. 
The spatial structure is a metapopulation of patches that can be of different habitat qualities. We label by $i=1,\ldots,K$ the patches so that the model can  be described by a labeled finite graph with weighted oriented edges. Vertices represent the patches, an oriented edge from vertex $i$ to vertex $j$ bears a weight $d_{ij}$ equal to the probability of dispersal from patch  $i$ to patch $j$. We let $D$ be the square matrix with generic element $d_{ij}$ and we call it \emph{dispersal matrix}. When $D$ has identical columns, we will speak of \emph{parent-independent migration}, or of a \emph{fully mixing} metapopulation.

We assume a simple asexual life cycle with discrete non-overlapping generations and no density-dependence. At each generation, as a net result of reproduction and survival (including survival to possible migration), all individuals of patch $i$, independently from one another, leave to the next generation a random number of individuals, called offspring, all distributed as some random variable $N_i$. The mean \emph{per capita} number of offspring in patch $i$ will be denoted by $m_i=\EE(N_i)$. In the second part of the paper, we will also enrich the model with a variable environment affecting simultaneously all patches. When the environment is in state $w$, we will denote by $m_i(w)$ the mean offspring in a patch of type $i$. 
   
Immediately after local growth, each individual from the new generation migrates independently, from patch $i$  to patch $j$ with probability $d_{ij}$. Since we assumed that mortality during dispersal is encompassed in the growth phase (see Discussion), we have  $\sum_{j=1}^K d_{ij}=1$ for all $i$, i.e., $D$ is a stochastic matrix. 

Reproduction, survival and dispersal probabilities are assumed not to depend on local densities. Thanks to this assumption of density-independence, and because the mean offspring numbers encompass migration-induced mortality, we have the following classification. If $m_i>1$, we say that patch $i$ is a \emph{source}, and if $m_i\leq 1$, we say that patch $i$ is a \emph{sink}.  

It will be convenient to assume that $m_1\geq m_2\geq \cdots\geq m_K$. In addition, the problem of persistence is  more interesting  in the case when $m_1>1\geq m_K$. Indeed, even in the presence of sinks, the metapopulation might persist thanks to local growth on sources replenishing sinks by dispersal. The case $m_K\leq  \ldots \leq m_1 \leq 1$ (resp.  $1<m_K\leq  \ldots \leq m_1$) leads trivially to extinction (resp. to persistence with positive probability).

\subsection{Two natural examples with one source type and one sink type} 

Let us describe two examples with two possible habitat qualities, one source type and one sink type, that will be treated as a special case throughout the paper:
$$
M:=m_1>1, \qquad m:=m_2=\cdots=m_K\leq 1.
$$ 

\begin{center}
\begin{figure}[ht]
\unitlength 1mm 
\linethickness{0.4pt}
\input{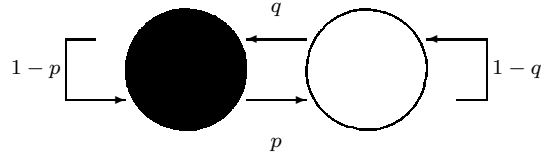}
\caption{Two patches of different qualities. The filled circle is a source and the empty circle is a sink. The arrow labels are the dispersal probabilities.}
\label{fig:two-patches}
\end{figure}
 \end{center}

 First, we will be interested in the simple case with two patches, patch $1$ with mean offspring  $M$, and patch $2$ with mean offspring  $m$. In this case, we will always use the simplified notation $p=d_{12}$ and $q=d_{21}$ (see Figure \ref{fig:two-patches}). 

\begin{center}
\begin{figure}[ht]
\unitlength 1.3mm 
\linethickness{0.4pt}
\input{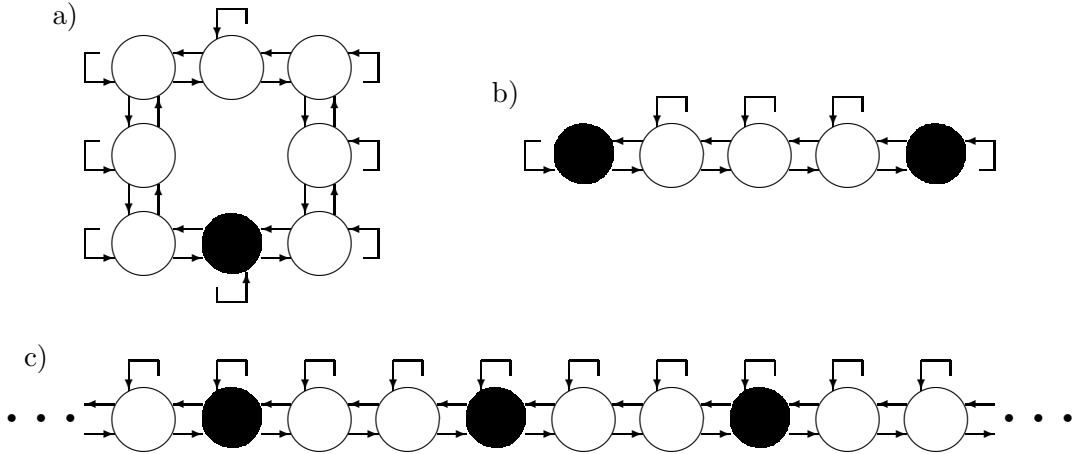}
\caption{Three examples where sources are only connected to sinks and two adjacent sources are connected by $n$ sinks; a) cyclic graph with one source, $n=7$; b) two sources, $n=3$; c) infinite periodic array, $n=2$, arrow labels (not represented) are also assumed periodic.}
\label{fig:arrays}
\end{figure}
 \end{center}

Second, we will consider the case  when each source is only connected to sinks and two adjacent sources are separated by an array of $n$ identical sinks. An example of such graph is the cyclic finite graph with one source and $n$ sinks, or two sources connected by $n$ sinks, or an infinite array with period $n$ (see Figure \ref{fig:arrays})...

\begin{table}[ht]%
\centering
\begin{tabular}{c|c}
 Notation & Interpretation  \\\hline
vertex $i$ & patch \\\hline
oriented edge  & probability of dispersal from patch $i$ to patch $j$ \\
 with weight $d_{ij}$ & \\ \hline
$m_i$ & mean number of offspring in a patch with habitat type $i$\\\hline
$M=m_1$ & mean growth rate in a source habitat (case $m_1>1\geq m_2=\cdots=m_K$)\\\hline
$m=m_2$ & mean growth rate in a sink habitat (case $m_1>1 \geq m_2=\cdots=m_K$)\\\hline
$p=d_{12}$ & probability of dispersal from the source to the sink (case of 2 patches) \\\hline
$q=d_{21}$ & probability of dispersal from the sink to the source (case of 2 patches)\\\hline
$X$ & random walk on the graph following the dispersal probabilities\\\hline
\end{tabular}
\caption{Notation.}
\end{table}

 \subsection{Method}
 
 The number of individuals located in  patch $i$  in generation $n$ is  denoted by $Z_n^{(i)}$. The process  $Z=(Z_n^{(i)}, \ i=1,\cdots,K, \   n \geq 0)$ is a \emph{multitype Galton--Watson process}. It is known from the mathematical literature \cite{AHBook, ANBook} that either  the population becomes extinct or it grows exponentially (under Assumption \textbf{(A1)} below). More specifically,  we see that $m_{i}d_{ij}$ is equal to the mean number of offspring of an individual living in patch $i$ which will land into patch $j$ in one time step, and therefore we call \emph{mean offspring matrix}  the matrix $A$ defined as
$$
A:=(m_{i}d_{ij} \  :  \ i,j=1,\cdots,K).
$$
The maximal eigenvalue (see e.g. \cite{SenBook}) of $A$ is the \emph{long-term growth rate}, or simply \emph{growth rate} of the metapopulation. Indeed \cite{AHBook, ANBook}, if $\rho\le 1$, the metapopulation dies out with probability $1$, and if $\rho>1$, the metapopulation can survive with positive probability, in which case 
$$
\frac{Z_n}{\rho^n}\stackrel{n\rightarrow \infty}{\longrightarrow} W ,
$$
where $W$ is a component-wise non-negative  and finite random vector (under Assumption \textbf{(A2')} below).

In the case of a fully mixing metapopulation, all column vectors of $D$ are identical to some vector $\delta$, say. 
Then $A$ is the rank 1 matrix $A= \delta \mu$, where $\mu$ is the row vector $\mu:=(m_1,\ldots, m_K)$. In this case, the spectral approach is straightforward, since $A^n = (\mu \delta)^{n-1} A$, so that $\rho=\mu\delta =\sum_{j=1}^K m_j\delta_j$, and $\delta$ and $\mu$ are respectively right and left eigenvectors of $A$ associated with $\rho$.

We call random disperser a single individual who moves on the graph at discrete time steps following the dispersal probabilities. In other words, if $X_n$ denotes the position of such a random disperser after $n$ time steps, then $(X_n)$ is the Markov chain with transition matrix $D$
$$
\PP(X_{n+1}=j\mid X_n=i) = d_{ij}.
$$
The goal of this paper is to display new persistence criteria, along with results regarding the asymptotic growth rate and the asymptotic fraction of time spent in each patch (by an individual taken at random in the surviving population). In contrast with the method involving the maximal eigenvalue of the mean offspring matrix, this one can yield quite simple, interpretable and partially explicit criteria. In addition, these criteria decouple the contributions of dispersal and demography on population survival. 
In a number of remarks, we will also provide sufficient conditions for survival which are explicit, in particular in the case of a fully mixing metapopulation. 

\begin{center}
\begin{figure}[ht]
\unitlength 1mm 
\linethickness{0.4pt}
\input{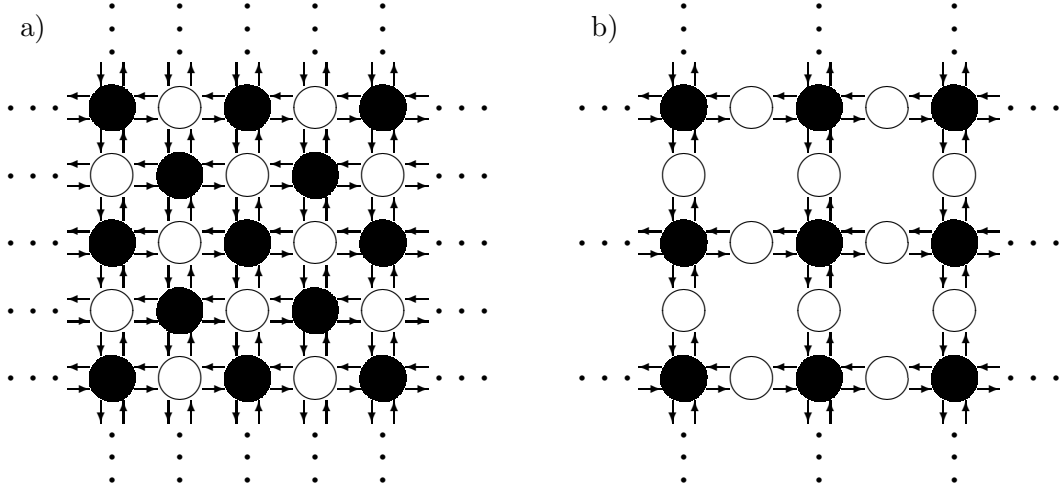}
\caption{Two examples of finite-transitive graphs   ; a) the chessboard; b) a square grid where four-degree vertices are sources separated by $n$ sinks (here $n=1$).}
\label{fig:transitive}
\end{figure}
\end{center}

\begin{center}
\begin{figure}[ht]
\unitlength 1mm 
\linethickness{0.4pt}
\input{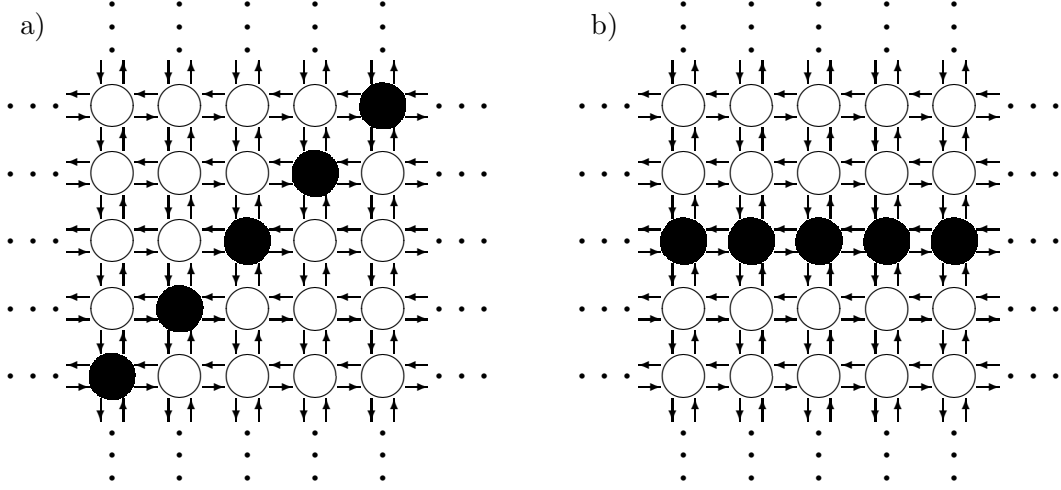}
\caption{Two examples of graphs embedded in $\ZZ^2$ which have an infinite motif (they still enjoy transitivity, but not finite) ; a) a diagonal of sources; b) a horizontal array of sources.}
\label{fig:other-transitive}
\end{figure}
\end{center}

This approach is still valid when the graph is an infinite  graph which can be reduced to a finite graph by transitivity. It is then called finitely transitive. The associated finite graph is called a motif, which is repeated to obtain the whole graph in such a way that the graph seen from any motif looks the same.
 A practical example is given by  sources with the same quality  connected by corridors of identical sinks and of the same length (see Figure \ref{fig:arrays} for an example). A finite-transitive graph could also be  an (infinite) chessboard where whites are  sinks and blacks are  sources, the square lattice $\ZZ^2$ where sources have coordinates of type $(n,n)$ (diagonal) or of type $(n,0)$ (horizontal array), and so on (see Figures \ref{fig:transitive} and \ref{fig:other-transitive})...
 
\subsection{Assumptions} 

Here, we list the assumptions we will use throughout the paper.

{\bf \emph{(A1)}} For convenience, we do not consider  the degenerate case when any individual in patch $1$ leaves one single offspring a.s. Thus, we assume that $ \mathbb P(N_1=1)<1$.

{\bf \emph{(A2)}} 
We assume that in each habitat $i$ the offspring number $N_i$ has a finite first moment, i.e., $\mathbb E (N_i)<\infty$. It ensures that the total size of the population has a finite first moment
in every generation. 

{\bf \emph{(A2')}} 
We assume that the offspring number $N_i$ of any individual living  in patch $i$ satisfies $\EE(N_i\log^+ N_i)<\infty$ (finite $N\log N$ moment).

{\bf \emph{(A3)}} 
To get the growth rate of the population, we will need to assume that $m_i \not= 0$ for all $i=1\ldots K$.

{\bf \emph{(A4)}} \emph{Irreducibility}. For any ordered pair $(i,j)$ of habitat patches, there is an integer $n$, such that $\PP(X_n=j\mid X_0=i)\not=0$. That is, the random disperser can go in finite time from any patch to any other patch by using  edges with positive weight. 

{\bf \emph{(A5)}} \emph{Aperiodicity}. For any patch $i$, the greatest common divisor (GCD) of the set of times $n$ such that $\PP(X_n=i\mid X_0=i)\not=0$, is equal to one.  
As a simple example, we mention the case when the graph is irreducible and at least one loop-edge has a positive weight, that is, there is at least one patch in which the probability of staying put is nonzero.  \\

When the graph is  both irreducible and aperiodic, we say that it is \emph{strongly irreducible} or \emph{primitive}. Then the Markov chain $(X_n;n \ge 0)$ is strongly  irreducible  and both the matrices $D$ (and $A$ if {\bf (A3)}  holds) are strongly irreducible, which is equivalent to the existence of  $n_0>0$ such that all the coefficient values of $M^{n_0}$ are positive \cite{SenBook}.

Assumption \textbf{(A2')}  ensures the convergence of $Z_n/\rho^n$ to a non degenerate r.v. $W$ which is non-negative on the survival event. This r.v. has only positive components if \textbf{(A3)} and \textbf{(A4)} are also in force.

\section{A first result on global persistence}
\label{secpersist}
\subsection{General case}

We now give a criterion for metapopulation persistence in terms of the random disperser $X$. 
For that purpose, we assume from now on that the random disperser starts  in patch $1$  ($X_0=1$) and we denote by $T$ the first return time of the random disperser into patch $1$,  
$$
T:=\min\{n\ge 1: X_{n}=1\}.
$$
\begin{thm}
\label{eqn : persistence crit2} We assume { \bf (A1, A2, A4)}.  Then the population persists with positive probability iff
$$
m_1\EE\left(\prod_{n=1}^{T-1} m_{X_n}\right) > 1.
$$
\end{thm}

In the case of a fully mixing metapopulation, $D$ has all its columns equal to some column vector $\delta$ and it is known that $\rho=\sum_{j=1}^K \delta_j  m_j$ (see previous section). It is easy to see that in this case, $(X_n)$ is a sequence of i.i.d. random variables whose common distribution is given by $\delta$, so that $T$ is geometrically distributed with success parameter $\delta_1$, and 
$$
m_1\EE\left(\prod_{n=1}^{T-1} m_{X_n}\right) = m_1 \sum_{n\ge 0} \delta_1 (1-\delta_1)^n\ \left(\frac{\sum_{j=2}^K \delta_jm_j}{1-\delta_1}\right)^n = \frac{\delta_1 m_1}{1-\sum_{j=2}^K \delta_j m_j}
$$
which is larger than 1 iff $\sum_{j=1}^K \delta_j  m_j$ is larger than 1. Thus, we recover the criterion obtained with the spectral approach.

\begin{rem}
\label{rem:occupation times}
Observe that the expression given in the theorem can also be expressed as $\prod_{n=1}^{T-1} m_{X_n}=\prod_{i=2}^K m_{i}^{S_{T-1}(i)},$
where $S_n(i):=\#\{ 1\leq k\leq n : X_k=i\}$ is the \emph{time spent in habitat $i$ by time $n$}, that is $S_{T-1}(i)$ is  the number of times the random disperser has visited patches of habitat type $i$ strictly before time  $T$. The advantage of this alternative formulation is that it carries over to models expressed in continuous time.
\end{rem}

Let us prove this result and then give more challenging applications.

\paragraph{Proof.}
Let $a$ (ancestor) be some  individual placed at time $0$ in  patch $1$. Define $Y_1$ as the number of offspring of $a$ staying put in patch $1$. Now for any integer $n\ge 2$, let $Y_n$ denote the number of descendants of $a$ at generation $n$ living in patch $1$ and whose ancestors  at generations $1, 2,\ldots, n-1$ have all lived \emph{outside} patch $1$. 
Then set
$$
Y:=\sum_{n\ge 1}Y_n,
$$
that can be seen as the total number of descendants of $a$ who live in patch $1$  for the first time in their lineage (except $a$).

In the theory of random trees, this set of individuals belonging to one of the $Y_n$ individuals for some $n$, is called a stopping line. It is known \cite{Cha86} that a stopping line enjoys the extended branching property, in the sense that all the subtrees descending from distinct elements of a stopping line are i.i.d. copies of the tree (conditional on their types). Then the total numbers of descendants of each of the individuals of this stopping line who live in patch $1$  for the first time in their lineage, are independent and all follow the same law as $Y$. In addition, any individual in the tree is either an ancestor or a descendant of some element of the stopping line.  Therefore, the total number of descendants of $a$ living in patch $1$  is finite iff the branching process with offspring number distributed as $Y$ is finite, which is equivalent to the a.s. extinction of $Y$.  The bottomline is that the population persists in habitat $1$  with positive probability iff $\EE(Y)>1$.  Indeed, we have excluded the critical case when $\PP(N_1=1)=1$.
But this local persistence in habitat $1$ is equivalent to the global persistence since the graph is  irreducible.

Let us then compute $\EE(Y)$ to conclude. We first note  that for every $i=1,\ldots,K$,  
$$\EE( Z_{n}^{(i)}  )=\sum_{j=1}^K \EE( Z_{n-1}^{(j)})m_{j}d_{ji},$$
where  $Z_{n}^{(i)}$ denotes the number of individuals located in patch $i$ at generation $n$. We prove easily by induction that the number of individuals $Y^{(i)}_n$ in patch $i$ at generation $n$ which have avoided patch $1$ at generations $k=1,2, \cdots, n-1$ satisfies
$$\EE(Y_n^{(i)})=\sum_{j=2}^K \EE(Y_{n-1}^{(j)})m_jd_{ji}=\sum_{j_1,\ldots, j_{n-1} \in \{2,\ldots K\}} d_{1j_1}d_{j_1j_2}\ldots d_{j_{n-2}j_{n-1}}d_{j_{n-1}i} m_1m_{j_1}\cdots m_{j_{n-2}}m_{j_{n-1}}.$$
As $Y_n=Y_n^{(1)}$, we get
\begin{eqnarray*}
\EE(Y_n)&=& \sum_{j_1,\ldots, j_{n-1} \in \{2,\ldots K\}} d_{1j_1}d_{j_1j_2}\ldots d_{j_{n-2}j_{n-1}}d_{j_{n-1}1} m_1m_{j_1}m_{j_2}\cdots m_{j_{n-2}}m_{j_{n-1}} \\
&=& m_1\EE(1_{ T=n}m_{X_1}m_{X_{2}}\cdots m_{X_{n-2}}m_{X_{n-1}}).
\end{eqnarray*}
Adding that  $Y=\sum_{n\geq 1} Y_n^{(1)}$, we have 
$$
\EE(Y)=   m_1\sum_{n \ge 1 } \EE(1_{ T=n}m_{X_1}\cdots m_{X_{n-1}}) =m_1\EE(m_{X_1}\cdots m_{X_{T-1}}).$$
This yields the result.\hfill $\Box$

\subsection{Case of two habitat types}
Let us focus now on the special case when there are $2$ habitat types and the source is  solely connected to sinks:
$$M:=m_1 >1, \qquad m:=m_2=\cdots=m_{K}<1.$$
We denote by 
$$p=\sum_{j=2}^K d_{1j}$$
 the probability of dispersing  for an individual living in patch $1$. The \emph{per capita} mean offspring number sent out from a source at each generation is $Mp$.
Let $\sigma$ be the time of \emph{first visit of a sink} by the random disperser
$$
\sigma:=\inf\{n\ge 0: X_n \ne 1\},
$$
so that $\sigma$ is a geometric random variable with success probability $p$.
Next, let $S$ denote the waiting time (after $\sigma$) before the random disperser visits a source (this source might or might \emph{not} be the initial source patch $X_0$)
$$
S:=\inf\{n\ge 0: X_{\sigma +n}=1\}.
$$
The duration $S$ can be seen as the \emph{time spent in sinks between two consecutive visits of sources}. By using the first transition of the random disperser, we get 
$$
\EE\left(\prod_{i=1}^{T-1} m_{X_n}\right)=1-p+p\ \EE\left(m^S\right),
$$
so that the previous theorem reads as follows. 
\begin{prop}
\label{thm initial} We assume { \bf (A1, A2, A4)}. 
Then the population persists with positive probability iff
\begin{equation}
\label{eqn : persistence crit}
M(1-p) + eMp >1,
\end{equation}
where $e$ is the depleting rate due to the sink habitat in the graph, defined as
$$
e:=\EE\left(m^S\right)= \sum_{k\ge 1}m^k \ \PP(S=k).
$$
\end{prop}

\begin{rem}
If the average time spent in sinks has
$$
\EE(S)<\frac{M-1}{Mp(1-m)},
$$
then the population persists with positive probability. Indeed, the mapping $f:x\mapsto \EE(x^S)$ is convex so 
$$
e=f(m)\ge 1+f'(1)(m-1)=1-(1-m)\EE(S)>1-\frac{M-1}{Mp}=\frac{1-M(1-p)}{Mp},
$$
which yields $eMp+M(1-p)>1$.
\end{rem}






Let us check, in the simple case  when there are only one source and one sink ($K=2$ vertices), that criterion \eqref{eqn : persistence crit} is equivalent to the condition that the maximal eigenvalue $\rho$ of $A$ exceeds unity.  Here the mean offspring matrix $A$ is
$$
A=
\left(
\begin{array}{cc}
M(1-p)& Mp\\
mq & m(1-q)	
\end{array}\right) .
$$ 
The characteristic polynomial $C$ of this square matrix is 
$$
C(x)=(M(1-p)-x)(m(1-q) -x)-Mmpq.
$$
Either $M(1-p)> 1$ and the population living in the source ensures the persistence. Or 
$M(1-p) \leq 1$ and the
quadratic polynomial is convex and has non negative  derivative at $1$. Thus, its leading eigenvalue is greater than 1 iff $C(1)<0$, which reads
$$
\frac{Mp}{1-M(1-p)}> \frac{1-m(1-q)}{mq} .
$$
We recover \eqref{eqn : persistence crit} since here $S$ is geometric with success probability $q$, which yields
$$
e=\sum_{k\ge 1} q(1-q)^{k-1}m^k=\frac{mq}{1-m(1-q)} .
$$
Notice that even in this simple case where $A$ is a $2\times 2$ matrix, the computation of the leading eigenvalue is cumbersome, and we have used a trick to explicitly specify the persistence criterion.

\subsection{Example with pipes of identical sinks}
\begin{center}
\begin{figure}[ht]
\unitlength 1.3mm 
\linethickness{0.4pt}
\input{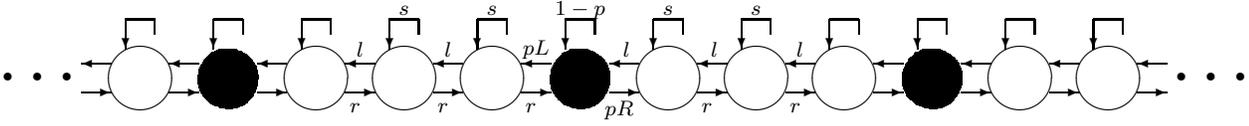}
\caption{A periodic pipeline where adjacent sources are separated by $n$ identical sinks ($n=3$).}
\label{fig:pipeline}
\end{figure}
\end{center}

Assume that the source  is a vertex of degree $2$ in the graph, connected to a \emph{left} sink and a  \emph{right} sink.  The probability of staying put on a source is still $1-p$, the probability of dispersing onto a left sink is $pL$, and the probability of dispersing onto a right sink is $pR$ (so that $L+R=1$). The sinks form a pipeline of $n$ adjacent sinks linking adjacent sources.  
The probability of staying put on a sink is always $s$,  the probability of dispersing from a sink onto one of its two neighboring sinks is $r$ in the left-to-right direction of the pipe, and $l$ in the right-to-left direction of the pipe (so that $q=l+r=1-s$). See Figure \ref{fig:cycle} for an example. 
This example will be directly  extended in the last Section to infinite  graphs, where pipelines of $n$ sinks periodically connect sources (see Figure \ref{fig:pipeline}). 

\begin{center}
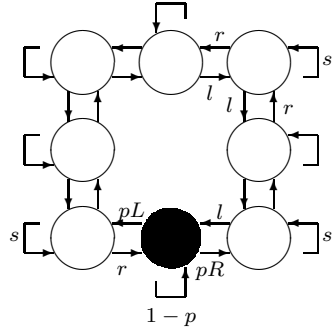
\begin{figure}[ht]
\unitlength 1.3mm 
\linethickness{0.4pt}
\begin{picture}(68,32.813)(-15,0)
\put(28,26.813){\circle{6}}
\put(37,26.813){\circle{6}}
\put(46,26.813){\circle{6}}
\put(28,17.813){\circle{6}}
\put(28,8.813){\circle{6}}
\put(46,17.813){\circle{6}}
\put(46,8.813){\circle{6}}
\put(34.8224,6.6354){\rule{4.3551\unitlength}{4.3551\unitlength}}
\multiput(35.7957,10.8781)(0,-4.893){2}{\rule{2.4086\unitlength}{.7628\unitlength}}
\multiput(36.3585,11.5284)(0,-5.714){2}{\rule{1.283\unitlength}{.2832\unitlength}}
\multiput(37.529,11.5284)(-1.4561,0){2}{\multiput(0,0)(0,-5.6425){2}{\rule{.3981\unitlength}{.2116\unitlength}}}
\multiput(38.0918,10.8781)(-2.8148,0){2}{\multiput(0,0)(0,-4.6157){2}{\rule{.6311\unitlength}{.4856\unitlength}}}
\multiput(38.0918,11.2512)(-2.5622,0){2}{\multiput(0,0)(0,-5.1402){2}{\rule{.3786\unitlength}{.2638\unitlength}}}
\multiput(38.0918,11.4025)(-2.4307,0){2}{\multiput(0,0)(0,-5.3577){2}{\rule{.2471\unitlength}{.1787\unitlength}}}
\multiput(38.358,11.2512)(-2.9565,0){2}{\multiput(0,0)(0,-5.0676){2}{\rule{.2406\unitlength}{.1912\unitlength}}}
\multiput(38.6105,10.8781)(-3.5699,0){2}{\multiput(0,0)(0,-4.4404){2}{\rule{.3489\unitlength}{.3102\unitlength}}}
\multiput(38.6105,11.0758)(-3.4538,0){2}{\multiput(0,0)(0,-4.7287){2}{\rule{.2329\unitlength}{.2031\unitlength}}}
\multiput(38.8469,10.8781)(-3.9179,0){2}{\multiput(0,0)(0,-4.3442){2}{\rule{.224\unitlength}{.214\unitlength}}}
\multiput(39.0651,7.6087)(-4.893,0){2}{\rule{.7628\unitlength}{2.4086\unitlength}}
\multiput(39.0651,9.9048)(-4.6157,0){2}{\multiput(0,0)(0,-2.8148){2}{\rule{.4856\unitlength}{.6311\unitlength}}}
\multiput(39.0651,10.4235)(-4.4404,0){2}{\multiput(0,0)(0,-3.5699){2}{\rule{.3102\unitlength}{.3489\unitlength}}}
\multiput(39.0651,10.6599)(-4.3442,0){2}{\multiput(0,0)(0,-3.9179){2}{\rule{.214\unitlength}{.224\unitlength}}}
\multiput(39.2628,10.4235)(-4.7287,0){2}{\multiput(0,0)(0,-3.4538){2}{\rule{.2031\unitlength}{.2329\unitlength}}}
\multiput(39.4382,9.9048)(-5.1402,0){2}{\multiput(0,0)(0,-2.5622){2}{\rule{.2638\unitlength}{.3786\unitlength}}}
\multiput(39.4382,10.171)(-5.0676,0){2}{\multiput(0,0)(0,-2.9565){2}{\rule{.1912\unitlength}{.2406\unitlength}}}
\multiput(39.5895,9.9048)(-5.3577,0){2}{\multiput(0,0)(0,-2.4307){2}{\rule{.1787\unitlength}{.2471\unitlength}}}
\multiput(39.7154,8.1715)(-5.714,0){2}{\rule{.2832\unitlength}{1.283\unitlength}}
\multiput(39.7154,9.342)(-5.6425,0){2}{\multiput(0,0)(0,-1.4561){2}{\rule{.2116\unitlength}{.3981\unitlength}}}
\put(40,8.813){\line(0,1){.3687}}
\put(39.977,9.182){\line(0,1){.1826}}
\put(39.949,9.364){\line(0,1){.1805}}
\put(39.909,9.545){\line(0,1){.1777}}
\put(39.859,9.722){\line(0,1){.1743}}
\put(39.797,9.897){\line(0,1){.1702}}
\put(39.725,10.067){\line(0,1){.1654}}
\put(39.643,10.232){\line(0,1){.16}}
\put(39.551,10.392){\line(0,1){.154}}
\put(39.449,10.546){\line(0,1){.1474}}
\multiput(39.337,10.694)(-.03005,.03508){4}{\line(0,1){.03508}}
\multiput(39.217,10.834)(-.03216,.03316){4}{\line(0,1){.03316}}
\multiput(39.088,10.967)(-.03414,.03112){4}{\line(-1,0){.03414}}
\put(38.952,11.091){\line(-1,0){.144}}
\put(38.808,11.207){\line(-1,0){.1508}}
\put(38.657,11.314){\line(-1,0){.1571}}
\put(38.5,11.411){\line(-1,0){.1628}}
\put(38.337,11.498){\line(-1,0){.1679}}
\put(38.169,11.576){\line(-1,0){.1723}}
\put(37.997,11.642){\line(-1,0){.1761}}
\put(37.821,11.698){\line(-1,0){.1792}}
\put(37.642,11.744){\line(-1,0){.1816}}
\put(37.46,11.777){\line(-1,0){.1834}}
\put(37.277,11.8){\line(-1,0){.1844}}
\put(37.092,11.812){\line(-1,0){.1848}}
\put(36.908,11.812){\line(-1,0){.1844}}
\put(36.723,11.8){\line(-1,0){.1834}}
\put(36.54,11.777){\line(-1,0){.1816}}
\put(36.358,11.744){\line(-1,0){.1792}}
\put(36.179,11.698){\line(-1,0){.1761}}
\put(36.003,11.642){\line(-1,0){.1723}}
\put(35.831,11.576){\line(-1,0){.1679}}
\put(35.663,11.498){\line(-1,0){.1628}}
\put(35.5,11.411){\line(-1,0){.1571}}
\put(35.343,11.314){\line(-1,0){.1508}}
\put(35.192,11.207){\line(-1,0){.144}}
\multiput(35.048,11.091)(-.03414,-.03112){4}{\line(-1,0){.03414}}
\multiput(34.912,10.967)(-.03216,-.03316){4}{\line(0,-1){.03316}}
\multiput(34.783,10.834)(-.03005,-.03508){4}{\line(0,-1){.03508}}
\put(34.663,10.694){\line(0,-1){.1474}}
\put(34.551,10.546){\line(0,-1){.154}}
\put(34.449,10.392){\line(0,-1){.16}}
\put(34.357,10.232){\line(0,-1){.1654}}
\put(34.275,10.067){\line(0,-1){.1702}}
\put(34.203,9.897){\line(0,-1){.1743}}
\put(34.141,9.722){\line(0,-1){.1777}}
\put(34.091,9.545){\line(0,-1){.1805}}
\put(34.051,9.364){\line(0,-1){.1826}}
\put(34.023,9.182){\line(0,-1){.9199}}
\put(34.051,8.262){\line(0,-1){.1805}}
\put(34.091,8.081){\line(0,-1){.1777}}
\put(34.141,7.904){\line(0,-1){.1743}}
\put(34.203,7.729){\line(0,-1){.1702}}
\put(34.275,7.559){\line(0,-1){.1654}}
\put(34.357,7.394){\line(0,-1){.16}}
\put(34.449,7.234){\line(0,-1){.154}}
\put(34.551,7.08){\line(0,-1){.1474}}
\multiput(34.663,6.932)(.03005,-.03508){4}{\line(0,-1){.03508}}
\multiput(34.783,6.792)(.03216,-.03316){4}{\line(0,-1){.03316}}
\multiput(34.912,6.659)(.03414,-.03112){4}{\line(1,0){.03414}}
\put(35.048,6.535){\line(1,0){.144}}
\put(35.192,6.419){\line(1,0){.1508}}
\put(35.343,6.312){\line(1,0){.1571}}
\put(35.5,6.215){\line(1,0){.1628}}
\put(35.663,6.128){\line(1,0){.1679}}
\put(35.831,6.05){\line(1,0){.1723}}
\put(36.003,5.984){\line(1,0){.1761}}
\put(36.179,5.928){\line(1,0){.1792}}
\put(36.358,5.882){\line(1,0){.1816}}
\put(36.54,5.849){\line(1,0){.1834}}
\put(36.723,5.826){\line(1,0){.1844}}
\put(36.908,5.814){\line(1,0){.1848}}
\put(37.092,5.814){\line(1,0){.1844}}
\put(37.277,5.826){\line(1,0){.1834}}
\put(37.46,5.849){\line(1,0){.1816}}
\put(37.642,5.882){\line(1,0){.1792}}
\put(37.821,5.928){\line(1,0){.1761}}
\put(37.997,5.984){\line(1,0){.1723}}
\put(38.169,6.05){\line(1,0){.1679}}
\put(38.337,6.128){\line(1,0){.1628}}
\put(38.5,6.215){\line(1,0){.1571}}
\put(38.657,6.312){\line(1,0){.1508}}
\put(38.808,6.419){\line(1,0){.144}}
\multiput(38.952,6.535)(.03414,.03112){4}{\line(1,0){.03414}}
\multiput(39.088,6.659)(.03216,.03316){4}{\line(0,1){.03316}}
\multiput(39.217,6.792)(.03005,.03508){4}{\line(0,1){.03508}}
\put(39.337,6.932){\line(0,1){.1474}}
\put(39.449,7.08){\line(0,1){.154}}
\put(39.551,7.234){\line(0,1){.16}}
\put(39.643,7.394){\line(0,1){.1654}}
\put(39.725,7.559){\line(0,1){.1702}}
\put(39.797,7.729){\line(0,1){.1743}}
\put(39.859,7.904){\line(0,1){.1777}}
\put(39.909,8.081){\line(0,1){.1805}}
\put(39.949,8.262){\line(0,1){.1826}}
\put(39.977,8.444){\line(0,1){.3687}}
\put(35.5,2.813){\line(1,0){3}}
\put(38.5,32.813){\line(-1,0){3}}
\put(52,16.313){\line(0,1){3}}
\put(52,7.313){\line(0,1){3}}
\put(52,25.313){\line(0,1){3}}
\put(22,19.313){\line(0,-1){3}}
\put(22,10.313){\line(0,-1){3}}
\put(22,28.313){\line(0,-1){3}}
\put(38.5,2.813){\vector(0,1){3}}
\put(35.5,32.813){\vector(0,-1){3}}
\put(52,19.313){\vector(-1,0){3}}
\put(52,10.313){\vector(-1,0){3}}
\put(52,28.313){\vector(-1,0){3}}
\put(22,16.313){\vector(1,0){3}}
\put(22,7.313){\vector(1,0){3}}
\put(22,25.313){\vector(1,0){3}}
\put(35.5,4.313){\line(0,-1){1.5}}
\put(38.5,31.313){\line(0,1){1.5}}
\put(50.5,16.313){\line(1,0){1.5}}
\put(50.5,7.313){\line(1,0){1.5}}
\put(50.5,25.313){\line(1,0){1.5}}
\put(23.5,19.313){\line(-1,0){1.5}}
\put(23.5,10.313){\line(-1,0){1.5}}
\put(23.5,28.313){\line(-1,0){1.5}}
\put(26.5,14.813){\vector(0,-1){3}}
\put(26.5,23.813){\vector(0,-1){3}}
\put(44.5,23.813){\vector(0,-1){3}}
\put(44.5,14.813){\vector(0,-1){3}}
\put(34,28.313){\vector(-1,0){3}}
\put(34,10.313){\vector(-1,0){3}}
\put(43,10.313){\vector(-1,0){3}}
\put(43,28.313){\vector(-1,0){3}}
\put(29.5,11.813){\vector(0,1){3}}
\put(29.5,20.813){\vector(0,1){3}}
\put(47.5,20.813){\vector(0,1){3}}
\put(47.5,11.813){\vector(0,1){3}}
\put(31,25.313){\vector(1,0){3}}
\put(31,7.313){\vector(1,0){3}}
\put(40,7.313){\vector(1,0){3}}
\put(40,25.313){\vector(1,0){3}}
\put(37,.75){\makebox(0,0)[cc]{\scriptsize $1-p$}}
\put(41,5.5){\makebox(0,0)[cc]{\scriptsize $pR$}}
\put(33,11.5){\makebox(0,0)[cc]{\scriptsize $pL$}}
\put(42,11.5){\makebox(0,0)[cc]{\scriptsize $l$}}
\put(32,5.5){\makebox(0,0)[cc]{\scriptsize $r$}}
\put(53,9){\makebox(0,0)[cc]{\scriptsize $s$}}
\put(21,9){\makebox(0,0)[cc]{\scriptsize $s$}}
\put(53,27){\makebox(0,0)[cc]{\scriptsize $s$}}
\put(42,29.5){\makebox(0,0)[cc]{\scriptsize $r$}}
\put(41,23.75){\makebox(0,0)[cc]{\scriptsize $l$}}
\put(49,22){\makebox(0,0)[cc]{\scriptsize $r$}}
\put(43,22.5){\makebox(0,0)[cc]{\scriptsize $l$}}
\end{picture}
\caption{A  pipeline where  $n$ identical sinks ($n=7$) connect the source to itself.}
\label{fig:cycle}
\end{figure}
\end{center}
We can compute exactly the depleting rate $e$ of the above kind of the corresponding infinite graph (finitely transitive). Let $\lambda> 1>\mu$ be the two ordered solutions to
$$
mrx^2-(1-ms)x+ml=0\qquad x\ge 0.
$$
Then
$$
\lambda\mu = \frac{l}{r} \qquad\mbox{ and }\qquad \lambda+\mu = \frac{1-ms}{mr} .
$$
\begin{prop} The depleting rate $e$ is equal to
$$
e=\frac{\lambda^{n}-\mu^{n}}{\lambda^{n+1}-\mu^{n+1}}(L+R\lambda\mu)
+
\frac{\lambda-\mu}{\lambda^{n+1}-\mu^{n+1}}(R+L(\lambda\mu)^n) .
$$
\end{prop}
\begin{rem}
In the two-patch case case ($n=1$), we recover 
$$
e=\frac{(1-s)m}{1-ms} .
$$
In the case when dispersal is isotropic $(l=r$), we get
$$
e=\frac{\lambda^{n}-\mu^{n}+\lambda-\mu}{\lambda^{n+1}-\mu^{n+1}} .
$$
Notice that in both previous cases, the depleting rate does not depend on $L$ or $R$.
In the case of one single source  and a large number of sinks ($n\tendinfty$), we get
$$
e=L\lambda^{-1} + R\mu.
$$
In the case of one single source and isotropic displacement, we then get $e=\lambda^{-1}=\mu$.
\end{rem}
\paragraph{Proof.} Consider a random walk $Y$ on $\{0,1,\ldots, n+1\}$, with displacement at each time step being $-1$ with probability $l$, $0$ with probability $s$, and $+1$ with probability $r$. Let $T_i$ denote the first hitting time of $i$ by $Y$, and set $T:=\min(T_0, T_{n+1})$ as well as
$$
a_k:=\EE(m^{T}\mid Y_0=k).
$$
Returning to the random disperser $X$ on the graph, it is easily seen that 
$$
\EE(m^S\mid X_\sigma \mbox{ is a left neighbour})= a_1,
$$
while
$$
\EE(m^S\mid X_\sigma \mbox{ is a right neighbour})= a_n,
$$
so that
$$
e=La_1+Ra_{n}. 
$$
Computations of $a_1$ and $a_n$ rely on the following recurrence relationship
$$
a_k = ms a_k + ml a_{k-1}+mr a_{k+1} \qquad k\in\{1,\ldots, n\},
$$
with boundary conditions $a_0= a_{n+1}=1$. This relation is obtained easily by considering the first transition of the walk $X$. \hfill $\Box$

\section{Growth rate and habitat occupation frequencies}

\subsection{General case}

If $\mathbf u$ denotes an individual in generation $n$, we define $H_k(\mathbf u)$ as the patch  occupied by the ancestor of $\mathbf u$ in generation $k\le n$. Then, for every $i\in \{1,\ldots,K\}$,
$$
F_i(\mathbf u):=\frac{1}{n}\,\#\{ 0 \leq k \leq n : H_k(\mathbf{u})=i\}
$$
is the \emph{occupancy frequency of patch $i$ by the ancestral line of $\mathbf{u}$}. We further denote by $\mathbf{U}_n$ an individual chosen randomly in the surviving population at generation $n$. We will see that the dispersal history of $\mathbf{U}_n$, as described by $(F_i(\mathbf{U}_n): i=1,\ldots,K)$ can be very different from that of a random disperser, since the ancestors of surviving individuals have better chance of having spent more time in sources than in sinks.

It is known that  the growth rate $\rho$ of the metapopulation is equal to  the maximal eigenvalue of $A$ (see \cite{AHBook, ANBook,  KLP97}). Moreover, the asymptotic occupancy frequencies of $\mathbf{U}_n$ are deterministic and can be expressed as  the product  of the right and left maximal eigenvectors associated to $\rho$.  We refer to \cite{Jag93, JN96, HRW02},  to Theorem 3.1 and 3.2 in  \cite{GB03} in continuous time and to  \cite{KLP97} in discrete time.
In this section, we want to give an alternative characterization to these quantities in terms of the random disperser and show an application.
 
To that purpose, we use the transition matrix $D$ of the random disperser $X$ on the graph. We denote by $F_i(X_n)$ the \emph{occupancy frequency of patch $i$ by the random disperser $X$ by time $n$}
$$
F_i(X_n):=\frac{1}{n}\,\#\{ 0 \leq k \leq n : X_k=i\}.
$$
By assumptions ({\bf A4, A5}), the random disperser has a stationary probability on $V$ (stochastic equilibrium), that we denote by $u=(u_i  : i=1,\ldots,K)$, which is the unique positive solution to $uD=u$. By the ergodic theorem, we also know that with probability 1,
$$
F_i(X_n)\stackrel{n\rightarrow \infty}{\longrightarrow} u_i .
$$
A typical single disperser will therefore occupy patch $i$ with asymptotic frequency $u_i$. This may not be the case of the ancestors of surviving individuals, whose paths must have favoured source patches. There is a trade-off between the pay-off in terms of fitness, gained by visiting source patches, and the cost in terms of likelihood, paid by deviating from the typical dispersal behavior. This trade-off is particularly obvious if we consider the case of a perfectly unsuitable habitat patch where the mean offspring is zero. In this case, the path followed by the ancestors of a surviving individual will necessarily have avoided this patch. Nevertheless, the asymptotic occupancy of this patch by a random disperser must be nonzero by the irreducibility assumption {\bf (A4)}.

There is a way of quantifying both the cost and pay-off of deviating from the typical dispersal behavior, that is, of having asymptotic occupancy frequencies $f=(f_i  : i=1,\ldots,K)$, where $f$ is a given element of the set $\cal F$ of non-negative frequencies on the graph
$$
\mathcal F:=\left\{f=(f_i :  i=1,\ldots,K) : \ f_i\geq 0, \ \sum_{i=1,\ldots,K} f_i=1\right\}.
$$

First, the probability that a random disperser has occupancy frequencies close to some given $f$ by time $n$ decreases exponentially with $n$ at rate $I(f)$, which can thus be interpreted as the cost of the $f$-occupancy scheme:
\begin{equation}
\label{def I}
I(f):=\sup\left\{\sum_{i=1,\ldots,K} f_i\log(v_i/(vD)_i) : v\gg0\right\},
\end{equation}
where $v\gg0$ denotes a positive row vector, that is, $v_i>0$ for each $i=1,\ldots,K$.
Indeed, large deviations theory \cite{DZBook, dHBook} ensures that for any $\epsilon\ll 1$,  as $n\to\infty$, we have
\begin{eqnarray}
\label{frqc1}
\PP(   f_i-\epsilon \leq F_i(X_n) \leq f_i+\epsilon \mbox{ for all } i=1, \cdots, K) \asymp \exp(-n I(f)).
\end{eqnarray}
We refer the reader to   Section \ref{proof} for a more  rigorous formulation. 
Taking $v=u$ in \eqref{def I} shows that $I(f)$ is of course always non-negative. This function is also convex.  When $f=u$, one can easily check that each partial derivative of $I$ is zero 
and  it can be proved that the supremum in \eqref{def I} is attained for $v=u$, so that $I(u)=0$. This was indeed expected, since $f=u$ is the natural occupancy scheme of the random disperser.

Second, the reproductive pay-off of $f$ can be defined as the fitness of a non-random disperser with given $f$-occupancy scheme, that is
$$
R(f):=\sum_{i=1,\ldots,K} f_i\log(m_{i}) .
$$
Indeed, the total size of a population of individuals  all adopting this dispersal behavior can be seen to grow like $$\prod_{i=1}^K m_i^{nf_i}=\exp(nR(f)).$$

Thus the cost (in terms of likelihood) for a population to follow some  occupancy scheme is quantified by $I$ and the reproductive pay off by $R$. The best  strategy (regarding the growth of the population) is to have an asymptotic occupancy frequency $\varphi$ which maximizes the difference $R-I$. If this optimum is positive, then the population  survives with positive probability. In addition, the ancestral line of a randomly chosen surviving individual will have visited patch $i$ with frequency $\varphi_i$. These results are stated below. The last assertion indicates that this optimal occupancy scheme $\varphi= (\varphi_i,i=1\ldots, K)$ is always different from the natural occupancy scheme $u$ of one single random disperser, except when all habitat types have the same quality.

\begin{thm} \label{rateped} We assume {\bf (A2', A3, A4, A5)}. \\
The growth rate $\rho$ of the metapopulation is given by
\begin{eqnarray*}
\log (\rho) &=&\max\left\{ R(f) - I(f): f\in \mathcal F\right\}. 
\end{eqnarray*}
In addition, if $\rho>0$, for any patch $i=1,\ldots,K$, conditional on the population being alive at time $n$, the occupancy frequency of patch $i$ by the ancestral line of a randomly chosen individual $\mathbf{U}_n$ in the surviving population at time $n$, converges to $\varphi_i$ in probability :
$$ 
F_i(\mathbf{U}_n)\stackrel{n\rightarrow \infty}{\longrightarrow } \varphi_i,
$$
where the frequency vector $\varphi=(\varphi_i) \in \mathcal F$  is uniquely  characterized by
$$
\log (\rho) = R(\varphi) - I(\varphi).
$$
The occupancy frequency $\varphi$ coincides with the stationary distribution $u$ of  $X$ (if and) only if
$$
m_1=m_2=\ldots =m_K.
$$
\end{thm}
In the same vein, we refer to Theorem 3.3 in \cite{GB03} for a description of the lineage of surviving individuals for multitype branching processes in continuous time.
The irreducibility assumption  {\bf (A4) } is required to use Sanov's large deviation theorem. The first result is actually a  consequence of the functional version of Sanov's theorem.
The aperiodicity assumption {\bf (A5)} and the exclusion of the degenerated cases by {\bf (A3)} are used for the two additional results.
The assumption {\bf (A2)} would be enough for the first part but {\bf (A2')} is required for the additional results. 

The proof is deferred to Section \ref{proof}.

\subsection{The fully mixing case} 

In the case of a fully mixing metapopulation, $D$ has all its columns equal to some column vector $\delta$ and we have already seen that $\rho=\sum_{j=1}^K \delta_j  m_j$. It is also easy to see that $\delta$ is the stable geographic distribution of the population. As a first observation, note that here the stable geographic distribution $\delta$ is also the stationary distribution of the random disperser. Then the last part of the previous theorem ensures that, except when all $m_j$'s are equal, the stable geographic distribution is different from the occupation distribution $\varphi$ of random long-lived lineages.

We now use the approach developed in the previous subsection for an alternative computation of $\rho$. We will also determine the occupation frequency $\varphi$ of ancestral lineages.
We first compute the functional $I$. Here, for any row vector $v\gg 0$,
$$
\sum_{i=1}^K f_i\log(v_i/(vD)_i) = \sum_{i=1}^K f_i\log(v_i/\vert v\vert \delta_i),
$$
where $\vert v \vert=\sum_{j=1}^K v_j$. Then,  differentiating this functional with respect to $v_j$ yields $(f_j/v_j) - 1/\vert v\vert$. As a result,
$$
I(f) = \sum_{i=1}^K f_i\log(f_i/\delta_i).
$$
Then differentiating $R-I$ with respect to $f_j$ ($j\not=1$ for example, and $f_1=1-\sum_{j=2}^K f_j$) yields $\log(\delta_j m_j) - 1- \log(f_j)$. The bottomline is 
$$
\varphi_j = \frac{\delta_jm_j}{\sum_{i=1}^K\delta_i m_i} .
$$
Plugging this as the argument of $R-I$ yields
$$
\log(\rho)= R(\varphi)-I(\varphi) = \log\left(\sum_{i=1}^K\delta_i m_i \right),
$$
which was the expected result. 

\section{Fluctuating  environments}

\subsection{General setting}

We now enrich our model with a fluctuating environment. The environment is embodied by a certain value $w$ which belongs to a finite set of states. We assume that the environment affects simultaneously all patches, but not necessarily  in the same way. We keep on assuming a simple asexual life cycle with discrete non-overlapping generations and no density-dependence. Now the environment is assumed to affect reproduction and survival, but not dispersal. Specifically, at each time step, conditional on the state $w$ of the environment, individuals reproduce independently according to some distribution which depends on the habitat type of their dwelling patch. We denote by $m_i(w)$ the mean offspring number of individuals dwelling in patch $i$ when the environment is in state $w$. 

Except in the last subsection, we will assume that the environment alternates periodically at each time step between two states (circadian cycle, seasons). Actually, the same method would allow to deal with any finite number of environmental states varying periodically.  

We call $e_1$ and $e_2$ the two possible states of the environment, so now we have $2K$  habitat qualities $m_i(e_j)$, for $i=1,\ldots, K$ and $j=1,2$. The Markov chain $Z=(Z_n^{(i)},   i=1,\ldots,K ;n\ge 0)$ is no longer time-homogeneous and is called \emph{multitype branching process in varying environment} \cite{HJVBook}. However, restricting the observation of the metapopulation to times when the environment is in the same state allows to adapt the arguments of the previous section. Indeed, $(Z_{2n}^{(i)}, i=1,\ldots,K;n\ge 0)$ is  a multitype branching process with mean offspring matrix $A$ with generic element
$$
a_{ij}=\sum_{k=1,\ldots,K} m_{i}(e_1)d_{ik}m_{k}(e_2)d_{kj}  \qquad i,j=1,\ldots,K.
$$
This amounts to changing the stopping line of the previous section, which was made of descendants returning to the ancestor patch for the first time in their lineage, for the stopping line of descendants returning to the ancestor patch for the first \emph{even} time in their lineage.
In the following subsection, we treat the case of two patches and determine the global persistence criterion. We then handle the general case using the random disperser.

\subsection{Example with two patches and two periodic environments}

For convenience, even if the environment is now variable, the two patches are still called respectively the source (patch 1) and the sink (patch 2). The mean number of offspring in the source are denoted by $M_1=m_1(e_1)$ and $M_2=m_1(e_2)$. In the sink, they are denoted by $m_1= m_2(e_1)$ and $m_2=m_2(e_2)$.

\begin{thm}
A necessary and sufficient condition for global persistence is
$$
M_1M_2(1-p)^2+(M_1m_2+m_1M_2)pq+m_1m_2(1-q)^2 > \min\big(2,1+ M_1M_2m_1m_2(1-p-q)^2\big).
$$
\end{thm}
\begin{rem}
It is easy to find examples where both patches are sinks on average but the metapopulation survives with positive probability thanks to dispersal. Indeed, each patch is a sink if (and only if) $M_1M_2\leq 1$ and $m_1m_2\leq 1$. Assuming for example that $p=q=1/2$ and $m_1=m_2=m$, the global survival  criterion becomes $M_1M_2+m(M_1+M_2)+m^2 > 4$, which holds as soon
 as $m(M_1+M_2) > 4$. 
\end{rem}

\paragraph{Proof.} The mean offspring matrix of $(Z_{2n};n\ge 0)$ is given by 
\begin{eqnarray*}
a(1,1)&=& M_1M_2(1-p)^2+M_1m_2pq \\
a(1,2)&=& M_1m_2p(1-q)+M_1M_2(1-p)p \\
a(2,1)&=&m_1M_2q(1-p)+m_1m_2(1-q)q \\
a(2,2)&=&m_1m_2(1-q)^2+m_1M_2qp.
\end{eqnarray*}
The maximum eigenvector of the matrix $A=(a(i,j) : 1\leq i,j\leq 2)$ is the largest root
 of the polynomial
 $$
 x^2-(a(1,1)+a(2,2))x+a(1,1)a(2,2)-a(1,2)a(2,1).
 $$
So it is less than $1$ iff
$$
a(1,1)+a(2,2)+\sqrt{(a(1,1)-a(2,2))^2+4a(1,2)a(2,1)}\leq 2
$$
Then the criterion for a.s. extinction of $(Z_{2n} : n \in \NN)$ is 
$$
a(1,1)+a(2,2)\leq 2 \quad \text{and} \quad (a(1,1)-a(2,2))^2+4a(1,2)a(2,1)\leq (2-a(1,1)-a(2,2))^2.
$$
The second inequality becomes
$a(1,1)+a(2,2)\leq 1+a(1,1)a(2,2)-a(1,2)a(2,1)$, which gives
$$
M_1M_2(1-p)^2+M_1m_2pq+m_1m_2(1-q)^2+m_1M_2qp\leq 1+ M_1M_2m_1m_2(1-p-q)^2.
$$
This completes the proof.\hfill $\Box$

\subsection{Global persistence for more than two patches}
Here, we extend the previous result to the case of a general, finite graph. We want to  state a global survival criterion which generalizes Theorem  \ref{eqn : persistence crit2} to periodic environments. Assume again that the random disperser starts at time $0$ in patch $1$  and set $T$ the first \emph{even} time when the random disperser goes back to habitat $1$
$$
T:=\min\{n\ge 1: X_{n}=1 \mbox{ and $n$ is even}\}.
$$
By a direct adaptation of the proof of Theorem \ref{eqn : persistence crit2} replacing  $Z_n$ with $Z_{2n}$, we get the following statement.
\begin{thm} We assume {\bf (A1)} holds for at least one environment,  {\bf (A2)} holds for both environments and {\bf (A4)} holds.
Then the population persists with positive probability iff
$$
m_1\EE\left(\prod_{n=1}^{T-1} m_{X_n}(w_n)\right) > 1,
$$
where the sequence $(w_n : n\geq 1)$ can take one of the two values $(e_1,e_2,e_1,\ldots)$ or $(e_2,e_1,e_2,\ldots)$, depending whether the initial environment is $e_1$ or $e_2$. 
\end{thm}

\subsection{Rate of growth and habitat occupation frequency }

The generalization to periodic environments of the results of the previous section can be achieved by changing the state-space $\{1,\ldots,K\}$ of the random disperser to the state-space of oriented edges of the graph, i.e., ordered pairs of vertices 
$$
 \mathcal{E}:= \{1,\ldots,K\}^2.
 $$ 
Denote by $B$ the transition matrix of the Markov chain $(X_{2n}, X_{2n+1};n\ge 0)$, which indeed takes values in $ \mathcal{E}$. Then denote by $\mathcal{F}$ the set of frequencies indexed by $\mathcal{E}$
$$
\mathcal{F}:=\big\{(f_E, E\in \mathcal{E}) : \ f_E\geq 0, \quad  \sum_{E\in \mathcal{E}} f_E=1\big\},
$$
and define the new cost function $I:\mathcal{F}\rightarrow \RR$ as
$$
I(f):= \sup\left\{\sum_{E\in \mathcal{E}} f_E\log(v_E/(vB)_E) : v\gg 0\right\},
$$
where $v$ denotes a non-negative vector indexed by $\mathcal{E}$, such that $v\gg0$, that is, $v_E>0$ for all $E\in \mathcal{E}$. Also define the new pay-off function  $R:\mathcal{F}\rightarrow \RR$ as
$$
R(f):=\sum_{E=(i,j) \in \mathcal{E}} f_E\log(m_{i}(e_1)m_{j}(e_2)).
$$
We can also provide an expression of  $I$ in terms of the entropy  function using Theorem 3.1.13 in \cite{DZBook}.
The generalization of Theorem \ref{rateped} can be stated as follows.
\begin{thm} \label{ratetwo} We assume that  {\bf (A2', A3)} hold for both environments and {\bf (A4, A5)} hold. \\
The growth rate $\rho$ of the metapopulation is given by
\begin{eqnarray*}
2\log (\rho) &=&\max\{ R(f) -I(f) : f\in \mathcal F\}.
\end{eqnarray*}
In addition, for any patch $i\in \{1,\ldots,K\}$, conditional on the population being alive at time $n$, the frequencies of occupation of patch $i$ by the ancestral line of a randomly chosen individual $\mathbf{U}_n$ in the surviving population at time $n$, converges  in probability :
$$
F_j(\mathbf{U}_n)\stackrel{n\rightarrow \infty}{\longrightarrow } \sum_{i \in \{1,\ldots,K\}} \varphi_{i,j}, 
$$
where the vector $(\varphi_{i,j} : (i,j)\in\mathcal{E})$ is characterized by
$$2\log (\rho) =  R(\varphi) -I(\varphi).$$
\end{thm}
 The proof follows that of Theorem  \ref{rateped}, with now 
$$\EE(\vert Z_{2n+1} \vert )= \EE\left(\prod_{k=0}^{n} m_{X_{2k}} m_{X_{2k+1}}\right)=\EE\left(\prod_{i=1}^K\prod_{j=1}^K m_i(e_1)^{S_n^{(1)}(i)}m_j(e_2)^{S_n^{(2)}(j)} \right),
$$
where $\vert Z_{n} \vert$ is the total number of individuals in source patches at generation $n$ and
$$
S_n^{(1)}(i)=\#\{ k\leq n : X_{2k}=i\}
,\quad S_n^{(2)}(i)=\#\{ k\leq n : X_{2k+1}=i\}
.
$$

\subsection{Example with fully  mixing patches}

We  extend the computations of the previous section to periodic environments. We focus on the fully mixing population : 
$$B_{(i,j)(k,l)}=d_{jk}d_{kl}=\delta_k\delta_l.$$

The mean matrix associated to the Galton Watson process $Z_{2n}$ is 
$$A_{ij}=m_i(e_1)\delta_j\left[\sum_{k=1}^K \delta_km_k(e_2)\right].$$
Thus the right and left eigenvectors  are still given by $\delta$  and $(m_1(e_1), \cdots m_K(e_1)$ and the spectral approach given previously can be followed readily. \\

Let us now focus on the approach given in the last Theorem.
As $(vB)_{E}=\vert v \vert \delta_{E_1}\delta_{E_2}$, the differentiation   of
$$\sum_{E\in \mathcal{E}} f_E\log(v_E/(vB)_E)=\sum_{E\in \mathcal{E}} f_E\log(v_E) - \log(\vert v_E\vert \delta_{E_1}\delta_{E_2})$$
with respect to $v_{E}$ yields the minimum. As in the previous section, we get
$$I(f)=\sum_{E\in \mathcal{E}} f_E\log\left(\frac{f_E}{\delta_{E_1}\delta_{E_2}}\right).$$
Then $R(f)-I(f)=\sum_{E\in \mathcal{E}} f_E\log\left(m_{E_1}(e_1)m_{E_2}(e_2)\delta_{E_1}\delta_{E_2}/f_E\right)$. We arbitrarily choose $E^0 \in \mathcal{E}$, so we can write  
$$f_{E^0}=1-\sum_{E \in \mathcal{E}_0} f_E $$ 
and make all other partial derivatives of $R-I$ equal $0$ when evaluated at $\varphi$. We get for every $E\ne E^0$ :
$$ \log\left(m_{E_1}(e_1)m_{E_2}(e_2)\delta_{E_1}\delta_{E_2}\right)+1-\log(\varphi_E)-\left[\log\left(m_{E^0_1}(e_1)m_{E_0^2}(e_2)\delta_{E^0_1}\delta_{E^0_2})\right)+1-\log(\varphi(E^0))\right]=0,$$
which gives the  habitat occupation frequencies
$$\varphi_E=\delta_{E_1}\delta_{E_2}m_{E_1}(e_1)m_{E_2}(e_2)\left[\sum_{E \in \mathcal{E}} \delta_{E_1}\delta_{E_2}m_{E_1}(e_1)m_{E_2}(e_2)\right]^{-1}.$$
We can also now deduce the growth rate $\rho$. 
We get 
$$ R(\varphi)-I(\varphi)=\log(\sum_{E \in \mathcal{E}} \delta_{E_1}\delta_{E_2}m_{E_1}(e_1)m_{E_2}(e_2))$$
and
$$\rho=\frac{1}{2}[ R(\varphi)-I(\varphi)] =\left[\sum_{i=1}^K \delta_i m_i(e_1)\right]^{1/2}\left[\sum_{i=1}^K \delta_i m_i(e_2)\right]^{1/2}.$$
This is the same growth rate as the one computed in \cite{JY98} for large populations (with two patches).


\subsection{Some comments on random environments}
A more natural way of modeling fluctuating environment in ecology is to assume random rather than periodic environment. The approach developed for periodic environments cannot be extended to random environments directly. Indeed, since the environment affects the whole metapopulation simultaneously, the randomness of environments correlates reproduction success in different patches. The process $(Z_n^{(i)}, i=1,\ldots,K;n\ge 0)$ counting the population sizes on each patch is now a \emph{multitype branching process in random environment} (MBPRE) \cite{AK71a, AK71b}.


Let us denote by $A(w)$ the mean offspring matrix (involving dispersal) in environment $w$. Specifically, the generic element $a_{ij}(w)$ of $A(w)$ is the mean offspring number of a typical individual dwelling in patch $i$ sent out to patch $j$ by dispersal, when the environment is $w$, so that
$$
A_{ij}(w)=m_{i}(w)\,d_{ij}.
$$
We will now assume that the state-space of environments is finite and that the sequence $(w_n: n \geq 0)$ of environment states through time is a \emph{stationary, ergodic sequence}, possibly autocorrelated, in the sense that the states need not be independent. 
Under this assumption, it is proved in \cite{FK60} (under the further assumption $\EE(\log^+\parallel A(w_0)\parallel )<\infty$, where expectation is taken w.r.t. the environment)  that the limit $\gamma$ of the sequence
$$ 
\frac{1}{n}\, \log\parallel  A(w_{n})A(w_{n-1})\ldots A(w_0)\parallel 
$$
exists with probability 1 and is deterministic, where $\parallel.\parallel$ denotes the maximum row sum of the matrix. This is interesting to us because it is further shown  in \cite{AK71a, Kap74, Tan77} (again under some further assumptions,  see Section \ref{classification}), that the extinction criterion and the growth rate of this MBPRE are respectively given by the sign and the value of $\gamma$, more specifically, $\gamma=\log \rho$. 

Unfortunately, this does not give a very explicit condition for global persistence. But again using the random disperser, we can give some sufficient conditions for survival. For simplicity, we turn our attention to the example of two patches and two environments $e_1$ and $e_2$. At any time step, the probability that the environment is in state $e_1$ is denoted by $\nu \in (0,1)$ (so the probability that the environment is in state $e_2$ is $1-\nu$). We show again that  the population may survive in sinks only.
As in the case of periodic environments, the mean number of offspring in the first patch (the source) is denoted by  $M_1=m_1(e_1)$ and $M_2=m_2(e_1)$. In the second patch (sink), they are denoted by $m_1= m_2(e_1)$ and $m_2=m_2(e_2)$. 


For the sake of simplicity, we state the results for the special case when the sequence is a Markov chain. We denote by $\alpha$ the transition from $e_1$ to $e_2$ and by $\beta$ the transition from $e_2$ to $e_1$. Then it is well-known that $\nu=\beta/(\alpha+\beta)$ is the asymptotic fraction of time spent in state $e_1$. The case of independent environments is recovered when $\alpha+\beta=1$. Note that as soon as $\alpha+\beta\not=1$, the sequence of environment states is auto-correlated.

\begin{prop} We have the following lower bound for the growth rate of the metapopulation.
$$
\log(\rho)\geq \nu\log(M_1)+(1-\nu) \log(m_{2})+ \nu\alpha\log(pq)+\nu(1-\alpha)\log(1-p)+(1-\nu)(1-\beta)\log(1-q).
$$
\end{prop}
\begin{rem}
Observe that this lower bound does not depend on $M_2$ and $m_1$. Again one can display examples where both patches are sinks but the metapopulation survives with positive probability in the presence of dispersal. Each patch is a sink if (and only if) $M_{1}^{\nu}M_{2}^{1-\nu}\le 1$ and $m_{1}^{\nu}m_{2}^{1-\nu}<1$. Actually one can manage to keep $\gamma >0$ while  $M_{1}^{\nu}M_{2}^{1-\nu}<1$, $m_{1}<1$ and $m_{2}<1$, for example with $M_{2}$ small and $M_1$ large for some fixed $p,q, m_1,m_2$. This corresponds to $e_2$ being a catastrophic environment in the source patch but the population survives  in patch $2$ when a  catastrophe occurs. 
\end{rem}


\paragraph{Proof.} We consider only  the subpopulation avoiding patch $1$ when the environment is equal to $e_2$. This means that this population reproduces with mean offspring number $M_1(1-p)$ in patch $1$ while the environment is  $e_1$. Each time the environment $e_2$ occurs, we consider the part of this population which has dispersed to patch $2$. This corresponds to a mean offspring number of $M_1p$. This population then stays in patch $2$ and  reproduces with mean offspring number $m_2(1-q)$ until the environment is again equal to $e_1$. We then consider the part of this population which goes back to patch $1$. This corresponds to a mean offspring number of $m_2q$.

Thus the patch of the ancestors of the individuals we keep is equal to $1$ (resp. $2$) if it lived  in environment $e_1$ (resp. $e_2$). Then at time $n$, the mean size of the population we consider is  equal to 
$$M_1^{N_1(n)}m_{2}^{N_2(n)}(1-p)^{N_{11}(n)}(1-q)^{N_{22}(n)}p^{N_{12}(n)}q^{N_{21}(n)}$$
where $N_i(n)$ ($i\in\{1,2\}$) is the number of times before generation $n$ when the environment is equal to $e_i$ 
 and $N_{ij}(n)$ ($i,j\in\{1,2\}$) is the number of one-step transitions of the environment from  $e_i$ to type $e_j$ until time $n$. By ergodicity, we know that  these quantities have  deterministic frequencies asymptotically. In the case of a Markovian sequence of environments, as $n\rightarrow \infty$,
$$
N_1(n)\sim \nu n, \quad N_2(n)\sim (1-\nu) n,
$$
and 
$$ 
N_{11}(n)\sim \nu(1-\alpha) n, \quad N_{12}(n)\sim \nu\alpha n,\quad N_{21}(n)\sim (1-\nu)\beta n, \quad N_{22}(n)=(1-\nu)(1-\beta)n.
$$
Using the growth of this  particular part of the whole population directly gives us a lower bound for $\gamma$:
$$
\gamma\geq \nu\log M_1+(1-\nu) \log m_{2}+ \nu\alpha\log p+(1-\nu)\beta\log q+\nu(1-\alpha)\log(1-p)+(1-\nu)(1-\beta)\log(1-q).
$$
Noticing that $(1-\nu)\beta=\nu \alpha$ completes the proof.  \hfill$\Box$

\begin{rem}
We could improve these results by considering more sophisticated strategies. For example, we could consider the subpopulation which stays in patch $1$ if (and only if) the number of consecutive catastrophes is less than $k$  and then optimize over $k$.

Actually, the proof relies on a stochastic coupling. Roughly speaking, the subpopulation we consider avoids the bad patches at the bad times and follows a (one type) branching process in random environment $e_{11}, e_{12}, e_{21}$ and $e_{22}$ respectively with stationary probabilities $\nu(1-\alpha)$, $\nu\alpha$, $(1-\nu)\beta$ and $(1-\nu)(1-\beta)$and mean offspring $M_1(1-p)$, $M_1p$, $m_2q$ and $m_2(1-q)$.  

Observe also that we can derive  a lower bound using the permanent of the mean matrix $M$ of the MBPRE from Proposition $2$ in \cite{BS09}. But this  lower   bound is not relevant for understanding the survival event in sinks only.
\end{rem}

\section{Metapopulation on infinite  graphs}
\label{transitivegraphs}

We now turn our attention to infinite graphs labeled by a countable set $F$. Each patch $P\in F$ has a type $i=\jmath(P) \in \mathbb N$  which gives its habitat quality, that is, the mean number of offspring in patch $P$ is equal to $m_{\jmath(P)}$. 

\begin{center}
\begin{figure}[ht]
\unitlength 1mm 
\linethickness{0.4pt}
\input{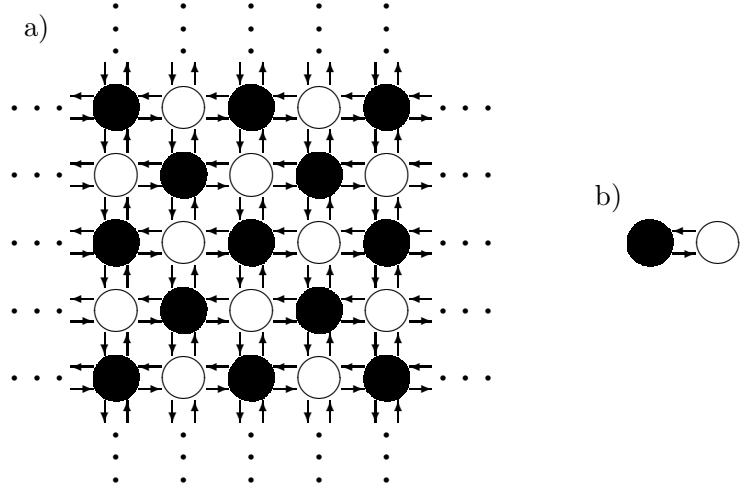}
\caption{a) The chessboard (with periodic arrow labels - not represented) is a  graph that can be collapsed into: b) a two-vertex graph (loop edges are not represented).}
\label{fig:chessboard}
\end{figure}
\end{center}

\begin{center}
\begin{figure}[ht]
\unitlength 1mm 
\linethickness{0.4pt}
\input{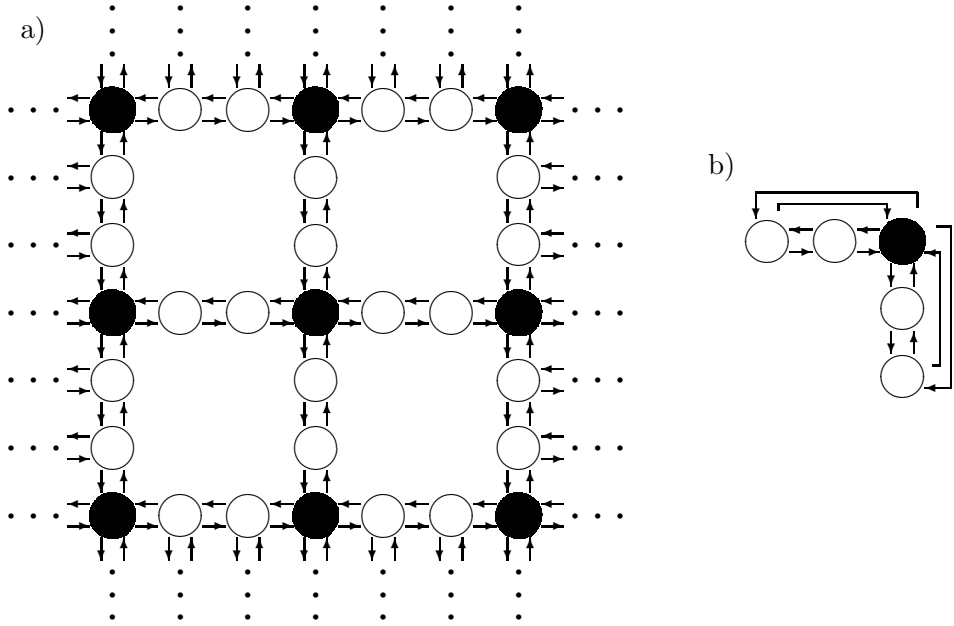}
\caption{ Star sources with $2d$ pipelines of sinks a) (and periodic arrow labels - not represented, here $d=2$) form a  graph  that can be collapsed into b) (a single source with $d$ pipelines, loop edges are not represented).}
\label{fig:stars}
\end{figure}
\end{center}

To generalize all the previous results, we require that the infinite graph has a \emph{finite motif}. Let us first provide the reader with some examples satisfying this assumption, before giving rigorous definitions. These examples are chosen among source-sink metapopulations with two habitat qualities (one source type and one sink type).
\begin{itemize}
\item An infinite linear periodic array of patches  (see Figure \ref{fig:arrays} for an example). The patches can then be labeled by $P\in F=\ZZ$ and the type of patch $P$ is equal to $\jmath(P)$, for some integers $N,K>0$ and a function  $\jmath : \ZZ \rightarrow [1,K]$ such that $\jmath(P+N)=\jmath(P)$ for every $P\in \ZZ$. The motif is a line of lenght $N$ with one source.
\item The chessboard (see Figure \ref{fig:chessboard}). The motif is composed of one source and one sink.
\item Star sources with $2d$ pipelines of sinks (see Figure \ref{fig:stars} for $d=2$). The motif is built by a source with $d$ pipelines of sinks of the same length.
\end{itemize}


Let us now specify mathematically these definitions. The oriented edges from $P$ to $Q$  are weighed by $d_{PQ}$. A  mapping $T$ of the graph is called  an \emph{isomorphism} if it conserves the types of the vertices as well as the weights of the oriented edges: $T$ is a bijection of $F$ such that  for all $P,Q \in F$,
$$ d_{T(P)T(Q)}=d_{PQ}, \qquad \jmath (T(P))=\jmath (P).$$
The associated equivalence relation $\sim$ between the patches of the graph is defined by 
$$P \sim P' \qquad \text{iff there exists an isomorphism } T \text{ of the graph such that } T(P)= P'.
$$
The class of a patch $P$ is defined as the equivalence class $Cl(P)=\{P' : P'\sim P\}$. Every patch of this class has the  type of  $P$. 
A graph for which there exists an isomorphism which is not the identity is called \emph{transitive}.

The collection of the distinct classes $(Cl(i) : i \in V)$ of a transitive graph form a partition of the patches of the graph. Such  subsets $V$ of patches are called  \emph{motifs}. 
 With a slight abuse of notation, the transition probabilities on  a motif $V$ are denoted by $(d_{PQ} : P\in V, Q\in V)$ where
$$d_{PQ}= \sum_{Q'\in Cl(Q)} d_{PQ'}$$
is constant in the same equivalence class. 
The initial graph can be seen as a family of copies of a motif properly connected. Observe that not all graphs with a regular structure are finitely transitive. For example, the cases illustrated by Figure \ref{fig:transitive} can be collapsed into a finite motif, but those given in Figure \ref{fig:other-transitive} (a diagonal or a horizontal array of sources in $\ZZ^2$) cannot. 
%
%
%
\begin{table}[ht]%
\centering
\begin{tabular}{c|c}
 Notation & Interpretation  \\\hline
vertex $P$ & patch \\\hline
set $V$ of vertices & motif \\\hline
type of $P$ $\jmath(P)=i$ & habitat quality of patch $P$ is $i$ \\\hline
oriented edge with weight $d_{PQ}$ & probability of dispersal from patch $P$ to patch $Q$ \\\hline
$m_i$ & mean number of offspring in a patch with habitat type $i$\\\hline
\end{tabular}
\caption{Notation.}
\end{table}
We can easily extend the results of  Section \ref{secpersist} to infinite graphs whose  motifs have  a finite number of sources. Thus, we consider the random disperser $X_n$ on the graph which follows the dispersal probabilities $(d_{PQ} : P,Q \in F^2)$.
To that purpose, we assume that the random disperser starts  in a patch of type $1$  ($\jmath(X_0)=1$) and we denote by $T$ the first return time of the random disperser into a (possibly different) patch of type $1$. 
$$
T:=\min\{n\ge 1: \jmath(X_{n})=1\}.
$$
\begin{thm} Assume that the metapopulation graph is transitive and that there is a motif with a finite number of sources. Assume {\bf (A1, A2, A4)} hold for this motif. 
Then the population persists with positive probability iff
$$
m_1\EE\left(\prod_{n=1}^{T-1} m_{\jmath(X_n)}\right) > 1.
$$
\end{thm}
We provide the idea of the proof. The number of individuals located in   patch $P$  in generation $n$ is  denoted by $N_n^{(P)}$. Collapsing the graph into some motif $V$, we denote by 
$$Z_n^{(P)}= \sum_{ P' \in Cl(P)} N_n^{(P')}$$
the total number of individuals in some habitat of the class of patch $P\in V$. 
Then $Z=(Z_n^{(P)}, \ P\in V;  n \geq 0)$ is a multitype Galton--Watson process (with possibly an infinite number of types).
We see that $m_{\jmath(P)}d_{PQ}$ is equal to the mean number of offspring of an individual living in patch $P$ which will land into patch $Q$ in one time step, and therefore we call mean offspring matrix  
$$
A:=(m_{\jmath(P)}d_{PQ} \  :  \ P,Q\in V),
$$
Then, everything  happens as if the metapopulation evolves on a  graph given by the motif. We follow  then the proof of Theorem \ref{fig:other-transitive}  to  prove that the criterion $m_1\EE(\prod_{n=1}^{T-1} m_{\jmath(X_n)}) >1$ is the persistence criterion in the  source habitat $1$. But using the facts that the number of sources in a motif is finite and that a motif is irreducible ensures that it is a global persistence criterion. \\

When the motif is finite, we get exactly the model described in the introduction and  derive the following counterpart of Theorem \ref{rateped}. We denote now by $D$ the transition matrix of the Markov chain $\jmath(X_n)$, and still use for $f=(f_P : P\in V)$
$$
I(f):=\sup\big\{\sum_{P\in V} f_P\log(v_P/(vD)_P) : v\gg0\big\},
$$
and 
$$
R(f):=\sum_{P \in V} f_P\log(m_{\jmath(P)}) .
$$
\begin{thm} Assume that the metapopulation graph is transitive with a finite motif. Assume {\bf (A1, A2', A3, A4, A5)} hold for this motif. 
 The growth rate $\rho$ of the metapopulation is given by
\begin{eqnarray*}
\log (\rho) &=&\max\left\{ R(f) - I(f): \sum_{P\in V} f_P=1, \ f_P \geq 0\right\}. 
\end{eqnarray*}
In addition, for any $P\in V$, conditional on the population being alive at time $n$, the occupancy frequency of a patch with type $\jmath(P)$ by the ancestral line of a randomly chosen individual $\mathbf{U}_n$ in the surviving population at time $n$, converges to $\varphi_i$ in probability
$$ 
F_i(\mathbf{U}_n)\stackrel{n\rightarrow \infty}{\longrightarrow } \varphi_i,
$$
where the vector $\varphi$ belongs to $\mathcal F$ and  is characterized by
$$
\log (\rho) = R(\varphi) - I(\varphi).
$$
The occupancy frequency $\varphi$ coincides with the stationary distribution $u$ of  $X$ only if $m_{\jmath(P)}=m_{\jmath(P')}$ for all patches $P,P'$.
\end{thm}

The results obtained for the periodic and random environments can be derived similarly for  transitive graphs.

\section{Discussion}

We have proposed two new techniques to derive simple criteria of persistence (by using the path of a random disperser) and to characterize the long-term growth rate of a single species in a source-sink metapopulation (by using large deviations for the path of the random disperser), as well as the occupancy frequencies of long-lived ancestral lineages. The expressions obtained thanks to these techniques decouple the contributions of reproduction/survival vs dispersal to the past, growth and persistence of these metapopulations. These techniques apply to a general class of stochastic, individual-based, stepping stone models of source-sink metapopulations, even when habitat quality is (not only variable in space but) variable through time, and even for (some) infinite metapopulations.

Growth rate and stable geographic distribution could  have as well been computed respectively as the maximal eigenvalue and associated eigenvector of the mean offspring matrix (encompassing dispersal). This well-known spectral approach can still more efficiently be used in general for the numerical computation of these quantities than the techniques presented here. However, the solution presented here has the advantage of giving a clear biological interpretation of the contributions of demography and dispersal to growth and stable distribution. In particular, we hope that this new presentation will help researchers in conservation biology to make informed decisions in planning  reintroduction strategies or in designing successfully protecting areas. Indeed, the benefit of our approach is to enable the ecologist to separate the effect on source-sink dynamics of changing the suitability of habitats (reproduction/survival scheme) or of changing the connectivity of the metapopulation (dispersal scheme).\\

Last, we want to indicate possible extensions of our ideas. 

First, we explain how migration-induced mortality has been encompassed in the growth phase, so that habitat patches $i$ with $m_i>1$ are sources in the strict sense of the word. Indeed, let $m_i'$ denote the mean number of offspring produced in patch $i$ \emph{before} dispersal, and $d_{ij}'$ denote the probability of migrating from $i$ to $j$ \emph{and} to survive to this migration event. Then $p_i=\sum_{i=1}^K d_{ij}'$ is the probability of survival to migration starting in patch $i$. It is easily seen that $m_i=p_im_i'$ is the mean offspring number including survival to migration and that  $d_{ij}=p_i^{-1} d_{ij}'$ are the dispersal probabilities of surviving individuals, so that considering only individuals surviving migration, everything happens as if growth with mean $m_i$ preceded dispersal with weights $d_{ij}$. These parameters $m_i$ and $d_{ij}$ are the ones we have used throughout the paper. Note that the matrix $A'$ with generic element $m_i' d_{ij}'$ is obviously equal to $A$. \\
On the other hand, we could as well have sticked to the framework where mortality-induced migration is not encompassed in the growth phase, in which case the dispersal matrix $D'$ is now sub-stochastic. Then the associated random disperser $X'$ would die with probability $1-p$ at each time step, and the criterion for persistence would have remained
$$
m_1'\EE\left(\prod_{n=1}^{T'-1} m_{X_n'}'\right) > 1,
$$
where $T'$ is the first return time to patch $1$ of the killed random walk $X'$, with the convention that the term inside the expectation is set to 0 when $T'$ is infinite (death of the walker before returning home).\\

Second, we can relax the assumption that dispersal behaviors of siblings are independent, provided there is no correlation between dispersal behaviors of different groups of siblings. Indeed, in this case, we just need to consider  the mean number $m_{ij}$ of offspring of an individual living in patch $i$ which go into patch $j$. Setting
    $$m_i=\sum_{j=1}^K m_{ij}, \qquad d_{ij}=m_{ij}/m_i$$
allows us to come back to our framework. In the case of finite multitype branching processes, the questions handled here are only linked to the mean offspring matrix. It makes such a procedure valid.

Third, if the model was expressed in continuous time then the disperser $X$ would be a time-continuous random walk, and the criterion for persistence to generalize is the criterion involving occupation times (see Remark \ref{rem:occupation times}).

Fourth, we mention that our ideas could also be adapted to age-structured populations, by considering the stopping line of descendants of a focal juvenile ancestor who are the first descending \emph{juveniles} to be born in the ancestor patch. Similarly, we can consider different (but a finite number) phenotypes or genotypes just by increasing the state space of the type of the branching process. Roughly speaking, habitat types are now replaced with a new composite type encompassing habitat, age,  phenotype, genotype... 


There is one further question that our methods could possibly solve. In infinite metapopulations, persistence can occur with the population failing to fill out the whole space, as in cases where dispersers always follow the same direction. This phenomenon is known as a dichotomy between local extinction and local exponential growth (conditional on global persistence). Our prediction of global persistence in infinite metapopulations relies on the study of the path of a random disperser started in a source patch, until the first time when it returns to a (possibly different) source patch. By making the difference between cases when the final patch is the same source patch or another source patch, we could display criteria for local persistence. Indeed, we know that there is \emph{global persistence} in a source-transitive patch iff
$$
m_1\EE\left(\prod_{n=1}^{T-1} m_{X_n}\right) > 1,
$$
where $T$ is the first return time to a source patch, but we conjecture that in this case there is \emph{local persistence} only if we also have
$$
m_1\EE\left(\prod_{n=1}^{\tau-1} m_{X_n}\right) > 1,
$$
where $\tau$ is the first return time to the very same source patch as initially.

\appendix

\section{Proof of Theorem \ref{rateped}}
\label{proof}
We assume {\bf (A2', A3, A4, A5)}. Thus the Markov chain $X$ is irreducible and evolves on a finite state space. 
The real number $I(f)$ is the cost for the habitat occupation frequencies associated with the random walk $X$ to equal $f$. It gives the geometric decrease
of the probability that the portion of time  spent in habitat $i$ until generation $n$ is close to $f_i$ :
\begin{eqnarray}
\label{frqc}
I(f_1,\ldots ,f_K) = \lim_{\epsilon \rightarrow 0} \lim_{n\rightarrow \infty} -\frac{1}{n} \log \PP(   f_k-\epsilon \leq F_k(X_n) \leq f_k+\epsilon).
\end{eqnarray}
 This result holds  when $I$ is finite. It is guaranteed by  Sanov's theorem (see e.g. \cite[Theorem 3.1.6 page 62]{DZBook}), which also ensures that  $I$   is convex  continuous.  The function $I$ is called the rate function associated to the path of the random walk $X$.\\

Finally, by {\bf (A2')}, the offspring distribution $N_i$ for an individual living  in patch $i$ satisfies $\EE(N_i\log^+ N_i)<\infty$, which ensures that for all $i,j=1,\ldots, K$,
$$\EE( Z_{1}^{j}\log^+  Z_1^{(j)} \ \vert \ Z_0^{(i)}=1, \ Z_0^{(k)}=0 \text{ for } k\ne i )<\infty .
$$
Using {\bf (A3, A4, A5)}, we have both $A$ and $D$ strongly irreducible. Then $Z_n/\rho^n$ converges to a non degenerate variable $W$, see e.g. \cite[Chapter 5, section 6, Theorem 1]{ANBook}. It is positive and finite on the survival event.

\paragraph{Expression of  the growth rate $\rho$ and  habitat occupation frequencies.}
 We start with one individual in patch $1$. 
Recalling that for every $i=1,\ldots ,K$,  $\EE( Z_{n}^{(i)}  )=\sum_{j=1}^K\EE( Z_{n-1}^{(j)})m_{j}d_{ji}$, we get by induction
$$\EE( Z_{n}^{(i)} )=\sum_{\substack{j_0=1, \ j_n=i, \\ \quad  1\leq j_1,\ldots ,j_{n-1} \leq K} } 
\prod_{k=0}^{n-1}m_{j_k}d_{j_k j_{k+1}}=\EE_1\big(1_{X_n=i}\prod_{i=0}^{n-1}m_{X_i}\big).$$
This yields
$$\EE(\vert Z_n \vert )= \EE\left(\prod_{i=0}^{n-1} m_{X_i}\right)$$
Denoting by $S_n(i)$ the number of  visits  of the random disperser $X$ in patch $i=1,\ldots V,$ :
$$S_n(i)=\#\{ k\leq n-1 : \jmath (X_k)=i\}=nF_i(X_n),$$
we deduce
\begin{eqnarray*}
\EE(\vert Z_n \vert )
&=&\EE\left(\prod_{i=1}^K m_{i}^{S_n(i)}\right) \\
&=&\int_{\mathcal{F}} \exp(n\sum_{i=1}^Kf_i\log(m_i))\PP(F_1(X_n) \in df_1,\ldots ,F_K(X_n)\in df_K). 
\end{eqnarray*}
Using $(\ref{frqc})$ and Laplace method, we get  
$$\log(\EE(\vert Z_n \vert)^{1/n}) \stackrel{n\rightarrow\infty}{\longrightarrow}  \max\{ \sum_{i=1}^K f_i\log(m_i) -I(f_1,\ldots ,f_K) : f\in \mathcal F\}.$$
This proves the first part of the result.\\

The  maximum of $h:=R-I$ is reached for a unique frequency $\varphi$, which means that there is a unique $\varphi$ such that $\log(\rho)=R(\varphi)-I(\varphi)$. This is the object of Proposition \ref{uniq}. Moreover  the partial derivatives of $h$ at $\varphi$ are zero. As $f_K=1-f_1-\ldots f_{n-1}$, for every $1\leq i\leq K-1$,
$$
\log(m_i)-\log(m_K)  -\frac{\partial}{\partial f_i} I_{\vert f=\varphi}=0.
$$
If $\varphi=p$, then the partial derivatives of $I$ at $\varphi$ are zero. This can be directly computed or deduced from   $I\geq 0$ and $I(p)=0$. This ensures that for every $1\leq i\leq K-1$, 
$\log(m_i)-\log(m_K)=0$, i.e. $m_1=m_2=\ldots =m_{K}$. This proves the third part of the theorem.

Finally, let us prove that the habitat occupation frequencies  of a typical individual are given by the vector $\varphi$. Let $\epsilon >0$ and $i \in\{1,\ldots,K\}$.
The individuals alive in generation $n$ are labelled by $\mathbf{u}_k$, $k=1,\ldots ,Z_n$.

Following the first part of the proof,
\begin{eqnarray*}
\EE(\sum_{k=1}^{Z_n} 1_{\vert F_i(\mathbf{u}_k)-\varphi_i \vert\geq \epsilon})&=& \EE\left(1_{\vert F_i(X_n)-\varphi_i \vert\geq \epsilon}\prod_{i=0}^n m_{X_i} \right)\\
&=&\EE\left(1_{\vert F_i(X_n)-\varphi_i \vert\geq \epsilon}\prod_{i=1}^K m_{i}^{S_n(i)}\right) \\
&=&\int_{\mathcal{F}} 1_{\vert f_i-\varphi_i \vert\geq \epsilon}\exp(n\sum_{i=1}^Kf_i\log(m_i))\PP(F_1(X_n) \in df_1,\ldots ,F_K(X_n)\in df_K).
\end{eqnarray*}
Using again  $(\ref{frqc})$ and the Laplace method, we get  
$$\frac{1}{n}\log \EE\big(\sum_{k=1}^{Z_n}
1_{\vert F_i(\mathbf{u}_k)-\varphi_i \vert\geq \epsilon}\big) \stackrel{n\rightarrow\infty}{\longrightarrow} C_{i,\epsilon},$$
with $C_{i,\epsilon}= \max\{ \sum_{i=1}^K f_i\log(m_i) -I(f_1,\ldots ,f_K) : f\in \mathcal F, \ \vert f_i-\varphi_i \vert\geq \epsilon\}.$
The  uniqueness of the argmax $\varphi$ ensures that 
$C_{i,\epsilon}< C_{i,0}$. Moreover the growth rate $C_{i,0}$ is equal to $\log(\rho)$ and
$$\frac{1}{n}\log \EE\big(\sum_{k=1}^{Z_n}
1_{\vert F_i(\mathbf{u}_k)-\varphi_i \vert\geq \epsilon}\big)-\frac{1}{n}\log \rho^n\stackrel{n\rightarrow\infty}{\longrightarrow} C_{i,\epsilon}-C_{i,0}<0.$$
Then 
$$\EE\big(  \frac{1}{\rho^n}.\sum_{k=1}^{Z_n} 1_{\vert F_i(\mathbf{u}_k)-\varphi_i \vert\geq \epsilon}\big)\stackrel{n\rightarrow\infty}{\longrightarrow}0.$$
In other words $\sum_{k=1}^{Z_n} 1_{\vert F_i(\mathbf{u}_k)-\varphi_i \vert\geq \epsilon}/\rho^n$ goes to $0$ in probability.
Adding that $Z_n \sim W \rho^n$ a.s. as $n\rightarrow \infty$ and $\{W>0\}=\{\forall n \in \NN, \  Z_n>0\}$ a.s. ensures that 
 $$ 1_{Z_n>0} \frac{1}{Z_n}.\sum_{k=1}^{Z_n} 1_{\vert F_i(\mathbf{u}_k)-\varphi_i \vert\geq \epsilon}\stackrel{n\rightarrow\infty}{\longrightarrow}0$$
in probability.  By dominated convergence,
$$ \EE\left(1_{Z_n>0} \frac{1}{Z_n}.\sum_{k=1}^{Z_n} 1_{\vert F_i(\mathbf{u}_k)-\varphi_i \vert\geq \epsilon} \right)\stackrel{n\rightarrow\infty}{\longrightarrow}0.$$
Then,  conditionally on $Z_n>0$, denoting  by $\mathbf{U}_n$ an individual chosen uniformly in generation $n$,
$$\PP\left(\vert F_i(\mathbf{U}_n)-\varphi_i \vert\geq \epsilon, \quad Z_n>0\right)\stackrel{n\rightarrow\infty}{\longrightarrow}0.$$
This proves that $ F_i(\mathbf{U}_n)\stackrel{n\rightarrow \infty}{\longrightarrow } \varphi_i$ in probability
and completes the proof.

\paragraph{Study of $I$ and uniqueness of argmax $R-I$.}
The supremum $I$ defined by
$$
I(f)=I(f_1,\ldots ,f_K):=\sup\{\sum_{j=1}^K f_j\log(u_j/(uD)_j) : u\in \RR^K, u\gg0\},$$
is reached for a unique unit positive vector. This means that there exists a unique $u(f)=(u_1,\ldots, u_K)$ such that
 $$I(f)=\sum_{j=1}^K f_j\log(u_j(f)/(u(f)D)_j), \quad u_1(f)+\cdots+u_K(f)=1, \quad  u_1(f)>0, \ldots, u_K(f)>0.$$ Indeed
this vector $u(f)$ realizes  a maximum  for $u\in \RR^K, u\gg0$ and thus satisfies for $j=1,\ldots ,K$,
\begin{equation} \label{idt}
 \frac{f_j}{u_j}-\sum_{i=1}^{K}d_{ji} \frac{f_i}{(uD)_i}=0.
\end{equation}
This equation  characterizes $u$, see  Exercise IV.9 page 46 in \cite{dHBook}, which ensures that $u(f)$ is uniquely defined. Note also that if $f_i$ is the stationary distribution, $u_j=f_j$ satisfies this equation since $(fD)_i=f_i$, so that $I=0$. \\

\begin{prop}\label{uniq}
There exists a unique $\varphi \in \mathcal{F}$ such that $\log(\rho)=R(\varphi)-I(\varphi)$.
\end{prop}
\paragraph{Proof.} 
We observe that $f\mapsto u(f)$ can be extended from $\mathcal{F}$ to $[0,\infty)^K\setminus\{0\}$ and can satisfy $(\ref{idt})$
on $[0,\infty)^K\setminus\{0\}$ by setting 
$$u(f)=u(f/ \parallel f \parallel ), \quad{where}   \ \ \parallel f \parallel = \sum_{i=1}^K f_i.$$ 
Then  $R-I$ can be extended   to  $[0,\infty)^K\setminus\{0\}$ with
$$(R-I)(f)= \sum_{j=1}^K \frac{f_j}{\sum_{k=1}^K f_k}\log(m_j(u(f)D)_j/u_j(f)),$$
and $(R-I)(\lambda f)=(R-I)(f)$ for every $\lambda \in (0,\infty)$.

Consider a vector $f$ which realizes the maximum of $R-I$ and does not belong to the boundary of $[0,\infty]^K$. Then the partial derivatives are zero 
 and for every $i=1, \ldots, K$,
$$\frac{1}{\sum_{k=1}^K f_k}.
\bigg[\log(m_i(uD)_i/u_i) +\sum_{j=1}^K f_j
\left(
\frac{\frac{\partial}{\partial f_i} (uD)_{j}}{(uD)_j}-\frac{\frac{\partial}{\partial f_i}u_j}{u_j}
\right)
\bigg]=\frac{\sum_{j=1}^K f_j\log(m_j(uD)_j/u_j)}{[\sum_{k=1}^K f_k]^2}.$$
Using (\ref{idt}) we get
$$
\sum_{j=1}^K f_j
\frac{\frac{\partial}{\partial f_i} (uD)_{j}}{(uD)_j}
=\sum_{j=1}^{K} \frac{f_j}{(uD)_j}\sum_{k=1}^K d_{kj} \frac{\partial u_k}{\partial f_i} 
= \sum_{k=1}^K   \frac{\partial u_k}{\partial f_i}  \sum_{j=1}^{K} \frac{f_jd_{kj}}{(uD)_j} 
=\sum_{k=1}^{K} \frac{\partial u_k}{\partial f_i}\frac{f_k}{u_k},
$$
so that
$$\log(m_i(uD)_i/u_i) =\frac{\sum_{j=1}^K f_j\log(m_j(uD)_j/u_j)}{\sum_{k=1}^K f_k}.$$
Observe that the right hand side does not depend on $i$, so that for every $i=1, \ldots, K$,
$$ (uD')_i=\alpha u_i,$$
where $D'_{ji}=d_{ji}m_i$ and $\alpha$ is a positive constant. Then $u$ is left eigenvector of $D'$ with positive entries.
Moreover $D'$ is strongly irreducible since $D$ is strongly irreducible  and $m_i>0$ for every $i$ by assumption (recall that $D$ is strongly irreducible if it is both irreducible and aperiodic, that is, if there exists $n_0\geq 1$ such that all the coefficients of $D^{n_0}$ are positive). 
Now Perron--Frobenius theory  ensures that there is a \emph{unique} left positive eigenvector  $u$
of $D'$ such that $\sum_{i=1}^K u_i=1$. Indeed, it is known \cite{HJBook} that two positive eigenvectors of a primitive matrix are colinear. This actually comes from the classical decomposition of $A^n$ using the maximum eigenvalue and the associated left and right eigenvectors.  Moreover, following the literature on large deviations \cite{dHBook}, (\ref{idt}) reads
$$
f_j=\sum_{i=1}^{K}D''_{ji} f_i,
$$
with $D''_{ji}=u_jd_{ji}/(uD)_i$.  Here again $D''$ is strongly irreducible since $D$ is strongly irreducible and both $u$ and $uD$ are positive vectors. Using again Perron--Frobenius theory guarantees the  uniqueness of the solution
$f$ such that $\sum_{i=1}^K f_i=1$. This ensures the uniqueness of the argmax of $R-I$ in the interior of $[0,\infty)^K$. We complete the proof by adding that there is at least one argmax in the interior of $[0,\infty)^K$ since we recall that the frequency occupation is  the product of the right and left eigenvectors of $A$ associate to $\rho$, which are both positive (using again Perron Frobenius theory with assumptions {\bf (A3, A4, A5)}).
If $\varphi_1$ and $\varphi_2$ realize the max of $R-I$, the concavity of this function ensures that so do all elements in the segment $[\varphi_1, \varphi_2]$, which is in contradiction with the uniqueness  in the interior of $[0,\infty)^K$ and completes the proof.
 \hfill$\Box$

\section{Classification Theorem for MBPRE}
\label{classification}
We consider here a multitype branching process in random environment $Z_n=(Z_n^{(i)} : i=1, \ldots,K)$ whose mean offspring matrix is denoted by $A=A(w)$, where $w=(w_0,w_1,\ldots)$ is the environment.\\
We introduce the extinction probability vector in environment $w$  starting from one individual in habitat $i$:
$$q_i(w)=\lim_{n\rightarrow \infty} \PP( \vert Z_n \vert =0 \ \vert \ w, \  Z_0^{(i)}=1, \  Z_0^{j}=0 \ \text{if } j\ne i).$$
\begin{prop}[\cite{Tan77}, Theorems 9.6 and 9.10] Assuming that
$$\PP( q_i(w)<1 : i =1,\ldots, K )= 1 \ \ \text{or} \ \ \PP( q_i(w) = 1 :  i =1,\ldots, K)= 1, \qquad (*)$$
we have   
\begin{itemize}
\item If $\gamma<0$, then  the probability of extinction is equal to $1$ for almost every $w$.
\item If $\gamma=1$, then
\begin{itemize}
\item either  for every $m \geq 1$, 
$w$-a.s., there exists  $1\leq i\leq K$ such that  

$\PP(\vert Z_m \vert > 1  \ \vert \ \ w, \ Z_0^{(i)}=1, \  Z_0^{(j)}=0 \ \text{if } j\ne i)=0$.
\item or  $q_i(w)=1$ $w$-a.s..
\end{itemize} 
\end{itemize}
Assuming that there exist integers $N,L>0$   such that $\PP( \forall 1\leq i,j\leq K, \ (A_N\cdots A_0)_{ij}\not=0   ) = 1$ and
$\vert \EE (\log (1-\PP(Z_K = 0 \ \vert \ Z_0^{(L)}=1) ) \vert <\infty$,  then $(*)$ is satisfied and
\begin{itemize}
\item If $\gamma>0$,   then $w$-a.s $q_i(w)<1$ for every $i=1,\cdots,K$, and 
$$\PP\big(\lim_{n\rightarrow \infty} n^{-1}\log(\vert Z_n\vert)=\gamma \ \vert  \ w,  \ Z_0^{(i)}=1, \  Z_0^{j}=0 \ \text{if } j\ne i \big)=1-q_i(w).$$
\end{itemize} 
\end{prop}
Moreover thanks to \cite[Theorem 9.11]{Tan77}, if all the coefficients of the matrix $A$ 
are positive and bounded, i.e.,
$$\exists 0<c,c'<\infty, \ c\leq \inf_{1\leq i,j\leq K} A_{i,j}\leq \sup_{1\leq i,j\leq K} A_{ij}\leq <c' \quad \text{a.s.},$$
then
$$Z_n= O(\| A_{n-1}\cdots A_0 \|)  \quad  \text{a.s.}$$

\paragraph{Acknowledgement.} We warmly thank Sebastian J. Schreiber and another (anonymous) reviewer for their feedback. This work was funded by project MANEGE `Mod\`eles
Al\'eatoires en \'Ecologie, G\'en\'etique et \'Evolution'
09-BLAN-0215 of ANR (French national research agency).

\bibliographystyle{abbrv}
\bibliography{myrefs} 



\end{document}